\documentclass[a4paper,12pt,reqno]{amsart}
\usepackage[a4paper,margin=1.9cm,top=2.5cm,bottom=2.5cm,centering,vcentering]{geometry}
\usepackage{hyperref,graphicx,tikz,amsmath,amssymb,amsfonts}

\definecolor{webgreen}{rgb}{0,.5,0}
\definecolor{webbrown}{rgb}{.6,0,0}
\definecolor{red}{rgb}{1,0,0}

\usepackage[normalem]{ulem}

\usepackage{color}

\usetikzlibrary{decorations.pathreplacing}

\newcommand{\mult}{\operatorname{mult}}

\renewcommand{\(}{\left\(}
\renewcommand{\)}{\right\)}
\renewcommand{\[}{\left\[}
\renewcommand{\]}{\right\]}

\numberwithin{equation}{section}
\theoremstyle{plain}
\newtheorem{theorem}{Theorem}[section]
\newtheorem{lemma}[theorem]{Lemma}

\newtheorem{remark}[theorem]{Remark}

\newtheorem{definition}[theorem]{Definition}

\newtheorem{problem}[theorem]{Problem}

\numberwithin{equation}{section}

\makeatletter
\def\proof{\@ifnextchar[{\@oproof}{\@nproof}}
\def\@oproof[#1][#2]{\trivlist\item[\hskip\labelsep\textit{#2 Proof of\
		#1.}~]\ignorespaces}
\def\@nproof{\trivlist\item[\hskip\labelsep\textit{Proof.}~]\ignorespaces}
\makeatother

\begin{document}
	\title{Parity biases in partitions and restricted partitions} 
%{\Large }

\author[K. Banerjee]{Koustav Banerjee}
\address{Research Institute for Symbolic Computation, Johannes Kepler University, Altenberger Strasse 69, A-4040 Linz, Austria}
\email{Koustav.Banerjee@risc.uni-linz.ac.at}

\author[S. Bhattacharjee]{Sreerupa Bhattacharjee}
\address{Department of Mathematics and Computer Science, University of Lethbridge, Lethbridge, AB Canada T1K 3M4} 
\email{bhatttacharjee.sreerupa@gmail.com}

\author[M. G. Dastidar]{Manosij Ghosh Dastidar}
\address{Technische Universit{\"a}t Wien, Wiedner Hauptstrasse 8–10/104, 1040 Wien, Austria} 
\email{gdmanosij@gmail.com}

\author[P. J. Mahanta]{Pankaj Jyoti Mahanta}
\address{Gonit Sora, Dhalpur, Assam 784165, India}
\email{pankaj@gonitsora.com}

\author[M. P. Saikia]{Manjil P. Saikia}
\address{School of Mathematics, Cardiff University, Cardiff, CF24 4AG, UK}
\email{manjil@saikia.in}



\date{\today.}

\maketitle

	\begin{abstract}
	Let $p_{o}(n)$ (resp. $p_{e}(n)$) denote the number of partitions of $n$ with more odd parts (resp. even parts) than even parts (resp. odd parts). Recently, Kim, Kim and Lovejoy proved that $p_{o}(n)>p_{e}(n)$ for all $n>2$ and conjectured that $d_{o}(n)>d_{e}(n)$ for all $n>19$ where $d_{o}(n)$ (resp. $d_{e}(n)$) denote the number of partitions into distinct parts having more odd parts (resp. even parts) than even parts (resp. odd parts). In this paper we provide combinatorial proofs for both the result and the conjecture of Kim, Kim and Lovejoy. In addition, we show that if we restrict the smallest part of the partition to be $2$, then the parity bias is reversed. That is, if $q_{o}(n)$ (resp. $q_{e}(n)$) denote the number of partitions of $n$ with more odd parts (resp. even parts) than even parts (resp. odd parts) where the smallest part is at least $2$, then we have $q_o(n)<q_e(n)$ for all $n>7$. We also look at some more parity biases in partitions with restricted parts.
\end{abstract}

\section{Introduction}\label{sec:intro}

In the theory of partitions, inequalities arising between two classes of partitions have a long tradition of study, for instance Alder’s conjecture \cite{alder} and the Ehrenpreis problem \cite{gea3} are the most famous examples in this direction. In recent years there have been a number of results about partition inequalities. For instance, work in this direction has been done by McLaughlin \cite{Mc}, Chern, Fu, and Tang \cite{Chern}, Berkovich and Uncu \cite{Berkovich} among others. Very recently, Kim, Kim, and Lovejoy \cite{kkl} have given interesting inequalities which show bias in parity of the partition functions. Further results on parity bias have been found by Kim and Kim \cite{kim_kim_2021}, and Chern \cite{Chern2}. Proofs of such results employ a wide range of techniques ranging from $q$-series methods, to combinatorial constructions and maps to classical asymptotic analysis. By parity bias we mean the tendency of partitions to have more parts of a particular parity than the other.

Let $p_{o}(n)$ (resp. $p_{e}(n)$) denote the number of partitions of $n$ with more odd parts (resp. even parts) than even parts (resp. odd parts). Kim, Kim, and Lovejoy \cite{kkl} proved that $p_{o}(n)>p_{e}(n)$ and conjectured that $d_{o}(n)>d_{e}(n)$ for all $n>19$ where $d_{o}(n)$ (resp. $d_{e}(n)$) denote the number of partitions into distinct parts having more odd parts (resp. even parts) than even parts (resp. odd parts). The primary goal of the present paper is to prove these two inequalities combinatorially.

In fact, our method can be amended to prove other results where biases in parity are found for restricted partitions. If $q_{o}(n)$ (resp. $q_{e}(n)$) denote the number of partitions of $n$ with more odd parts (resp. even parts) than even parts (resp. odd parts) where the smallest part is at least $2$, then we have $q_o(n)<q_e(n)$ for $n>7$ (see Theorem \ref{thm:reverse} below). These parity biases seems to also occur for more restricted partition functions and we also explore some of these themes towards the end with a few conjectures.

We define a partition $\lambda$ of a non-negative integer $n$ to be an integer sequence $(\lambda_{1},\dots,\lambda_{\ell})$ such that $\lambda_{1} \geq \lambda_{2} \geq \dots \geq \lambda_{\ell} >0$. We say that $\lambda$ is a partition of $n$, denoted by $\lambda \vdash n$ and $\sum_{i=1}^{\ell} \lambda_{i} = n$. The set of partition of $n$ is denoted by $P(n)$ and $|P(n)| = p(n)$. For $\lambda \vdash n$, we define $a(\lambda)$ to be the largest part of $\lambda$, $\ell(\lambda)$ to be the total number of parts of $\lambda$ and $\mult_{\lambda}(\lambda_{i}):=m_i$ to be the multiplicity of the part $\lambda_{i}$ in $\lambda$. We also use $\lambda= (\lambda_{1}^{m_{1}}\dots\lambda_{\ell}^{m_{\ell}})$ as an alternative notation for partition. For $\lambda \vdash n$ with $\lambda = (\lambda_{1},\dots,\lambda_{\ell})$ and $\mu \vdash m$ with $\mu = (\mu_{1},\dots,\mu_{\ell^{'}})$, define the union $\lambda \cup \mu  \vdash m+n$ to be the partition with parts $\{\lambda_{i},\mu_{j}\}$ arranged in non-increasing order. For a partition $\lambda \vdash n$, we split $\lambda$ into $\lambda_{e}$ and $\lambda_{o}$
respectively into even and odd parts; i.e., $\lambda = \lambda_{e} \cup \lambda_{o}$. We denote by $\ell_{e}(\lambda)$ (resp. $\ell_{o}(\lambda)$) to be the number of even parts (resp. odd parts) of $\lambda$ and $\ell(\lambda) = \ell_{e}(\lambda)+\ell_{o}(\lambda)$.

The following sets of partitions are of interest in this paper.
\begin{definition}
	\begin{equation*}
	\begin{split}
	D(n) & := \{\lambda \in P(n):\mult_{\lambda}(\lambda_{i})=1  \ \text{for all} \  i\},\\
	P_{e}(n) & := \{\lambda \in P(n):  \ell_{e}(\lambda)> \ell_{o}(\lambda)\},\\
	P_{o}(n) & := \{\lambda \in P(n):  \ell_{o}(\lambda)> \ell_{e}(\lambda)\},\\
	D_{e}(n) & :=  P_{e}(n) \cap D(n),\\
	D_{o}(n) & := P_{o}(n) \cap D(n),\\
	Q(n) & :=\{\lambda \in P(n): \lambda_i\neq 1 \ \text{for all} \ i\},\\
	Q_e(n) & := \{\lambda \in Q(n): \ell_e(\lambda)>\ell_o(\lambda)\},\\
	Q_o(n) & := \{\lambda \in Q(n): \ell_o(\lambda)>\ell_e(\lambda)\},\\
	DQ_{e}(n) & := Q_e(n) \cap D(n),\\
\text{and}\ \ DQ_{o}(n) & := Q_o(n) \cap D(n).
	\end{split}
	\end{equation*}
\end{definition}

\begin{definition}
	For all the sets defined above, their cardinalities will be denoted by the lower case letters. For instance, $|P_e(n)|=p_e(n)$, $|DQ_e(n)|=dq_e(n)$ and so on.
\end{definition}

Now, we state formally the main results proved in this paper.

\begin{theorem}[Theorem 1, \cite{kkl}]\label{thm1} For all positive integers $n \neq 2$, we have
	$$p_{o}(n)>p_{e}(n).$$
\end{theorem}

\begin{theorem}[Conjectured, \cite{kkl}]\label{thm:distinct}
	For all positive integers $n > 19$, we have
	$$d_{o}(n)>d_{e}(n).$$	
\end{theorem}

\begin{theorem}\label{thm:reverse}
	For all positive integers $n>7$, we have
	\begin{equation*}
	q_o(n)<q_e(n).
	\end{equation*}
\end{theorem}

For a nonempty set $S \subsetneq \mathbb{Z}_{\geq 0}$, define
\begin{equation*}
\begin{split}
P_e^S(n)&:=\{\lambda \in P_e(n): \lambda_i \notin S\}\\
\text{and}\ \ P_o^S(n)&:=\{\lambda \in P_o(n): \lambda_i \notin S\}.\\
\end{split}
\end{equation*}
Consequently, denote the number of partitions in $P_e^S(n)$ (resp. $P_o^S(n)$) by $p_e^S(n)$ (resp. $p_o^S(n)$).
The above definition leads us to the following results that describes not only the parity of parts but also its arithmetic by putting a constrain on its support.

\begin{theorem}\label{thm:missingset1}
	For all $n\geq1$ we have $$p_o^{\{2\}}(n)>p_e^{\{2\}}(n).$$
\end{theorem}

\begin{theorem}\label{thm:missingset2}
	If $S=\{1,2\}$, then for all integers $n>8$, we have
	$$p^{S}_o(n)>p^{S}_e(n).$$
\end{theorem}

Before we move on further, let us describe the fundamental principle behind proofs of Theorems \ref{thm1}-\ref{thm:missingset2}. Let $X$ and $Y$ be two given sets and our goal is to prove that $|Y| > |X|$. We choose a subset $X_{0}\ (\varsubsetneq X)$ and define an injective map $f: X_{0} \rightarrow Y$. Then to prove $|Y| > |X|$, it is enough to prove for a suitable subset $Y_{0}\ \varsubsetneq Y\setminus f(X_{0})$ with $|Y_{0}| > |X\setminus X_{0}|$. Throughout this paper, we follow the notation $x\mapsto y$ instead of writing $f(x)=y$ when the map $f$ is understood from the context.

The rest of the paper is organized as follows: in Section \ref{sec:proof-kkl} we give a combinatorial proof of the result of Kim, Kim and Lovejoy \cite{kkl}, in Section \ref{sec:distinct} we give a proof of the conjecture of Kim, Kim and Lovejoy \cite{kkl}, in Section \ref{sec:reverse} we prove reverse parity bias as stated in Theorem \ref{thm:reverse}, Section \ref{sec:missingset} presents the proofs of Theorems \ref{thm:missingset1} and \ref{thm:missingset2}, and finally in Section \ref{sec:final} we present a very short discussion drawing on Section \ref{sec:missingset}, by proposing further problems. The proofs of two preliminary lemmas (cf. Lemmas \ref{lem1} and \ref{lem2}) are given in Appendix \ref{sec:appendix}.

\section{Proof of Theorem \ref{thm1}}\label{sec:proof-kkl}
We begin by presenting the following two lemmas, used later in the proof of Theorem \ref{thm1}. For proofs, we refer to Appendix \ref{sec:appendix}.
\begin{lemma}\label{lem1}
	For all even positive integer $n$ with $n \geq 14$, we have
	\begin{equation*}
	\sum_{k=1}^{\frac{n-6}{2}}\Bigl\lfloor\frac{n-2k-2}{4}\Bigr\rfloor > 1 + \sum_{k=1}^{\lfloor \frac{n-2}{6}\rfloor}\Bigl\lfloor \frac{n-6k+2}{4}\Bigr\rfloor + \sum_{k=1}^{\lfloor \frac{n-6}{6}\rfloor}\Bigl\lfloor \frac{n-6k-2}{4}\Bigr\rfloor.
	\end{equation*}
\end{lemma}
\begin{lemma}\label{lem2}
	For all odd positive integer $n$ with $n \geq 11$, we have
	\begin{equation*}
	\sum_{k=1}^{\frac{n-5}{2}}\Bigl\lfloor\frac{n-2k-1}{4}\Bigr\rfloor > 1 + \Bigl\lfloor \frac{n-4}{4}\Bigr\rfloor+\sum_{k=1}^{\lfloor \frac{n-5}{6}\rfloor}\Bigl\lfloor \frac{n-6k-1}{4}\Bigr\rfloor + \sum_{k=1}^{\lfloor \frac{n-9}{6}\rfloor}\Bigl\lfloor \frac{n-6k-5}{4}\Bigr\rfloor.
	\end{equation*}
\end{lemma}

Let
\begin{equation*} 
\begin{split}
G^{0}_{e}(n) & := \{\lambda \in P_{e}(n):  \ell_{e}(\lambda)- \ell_{o}(\lambda)=1 \ \text{and} \ a(\lambda) \ \equiv 0 \ ( \text{mod} \ 2) \},\\
\overline{G^{0}_{e}}(n) & := \{\lambda \in G^{0}_{e}(n) : \lambda_{3} \geq 3\},\\
G^{1}_{e}(n) & := \{\lambda \in P_{e}(n):  \ell_{e}(\lambda)- \ell_{o}(\lambda)=1 \ \text{and} \ a(\lambda) \ \equiv 1 \ ( \text{mod} \ 2) \},\\
G^{2}_{e}(n) & := \{\lambda \in P_{e}(n): \ell_{e}(\lambda)- \ell_{o}(\lambda) \geq 2 \},\\
\text{and} \ \ G_{e}(n) & := G^{1}_{e}(n) \ \cup \ G^{2}_{e}(n).
\end{split}
\end{equation*}

We split the set $G_{e}(n)$ into the parity of length of partition as $G_{e}(n) = G_{e,0}(n) \ \cup \ G_{e,1}(n)$ with $G_{e,0}(n) = \{\lambda \in G_{e}(n) : \ell(\lambda) \equiv 0 \ ( \text{mod} \ 2)\}$, $G_{e,1}(n) = \{\lambda \in G_{e}(n) : \ell(\lambda) \equiv 1 \ ( \text{mod} \ 2)\}$ and let $\overline{G_{e}}(n) := G_{e,0}(n) \ \cup \ G_{e,1}(n) \ \cup \overline{G^{0}_{e}}(n)$. Therefore,
\begin{equation}\label{eq0}
P_{e}(n) \setminus \overline{G_{e}}(n) = \{\lambda \in G^{0}_{e}(n) : 0 \leq \lambda_{3} \leq 2 \}.
\end{equation}

We construct a map $f : \overline{G_{e}}(n) \rightarrow P_{o}(n)$ by defining maps $f|_{G_{e,0}(n)} = f_{1}$, $f|_{G_{e,1}(n)} = f_{2}$ and $f|_{\overline{G^{0}_{e}}(n)} = f_{3}$ such that $\{f_{i}\}_{1 \leq i \leq 3}$ are injective with the following properties
\begin{itemize}
	\item  $ f_{1}(G_{e,0}(n)) \ \cap \ f_{2}(G_{e,1}(n)) = \emptyset$,
	\item $ f_{1}(G_{e,0}(n)) \ \cap \ f_{3}(\overline{G^{0}_{e}}(n)) = \emptyset$, and
	\item $ f_{2}(G_{e,1}(n)) \ \cap \ f_{3}(\overline{G^{0}_{e}}(n)) = \emptyset$,
\end{itemize}
so as to conclude the map $f$ is injective. Then we will choose a subset $\overline{P_{o}}(n) \ \varsubsetneq P_{o}(n) \setminus f(\overline{G_{e}}(n))$ with $|\overline{P_{o}}(n)|>|P_{e}(n) \setminus \overline{G_{e}}(n)|$.

Let $\lambda \in G_{e,0}(n)$ with $\lambda_{e} = (\lambda_{e_{1}},\dots,\lambda_{e_{k}})$  and $\lambda_{o} = (\lambda_{o_{1}},\dots,\lambda_{o_{m}})$ where $k+m = \ell(\lambda)$. Since $\lambda \in G_{e,0}(n)$, $\ell(\lambda) =2r$ for some $r \in \mathbb{Z}_{>0}$ and $k >r$ because, $k-m \geq 1 \ \text{implies} \ 2k \geq k+m+1 = 2r+1$.

We define $f_{1} : G_{e,0}(n) \rightarrow P_{o}(n)$ by $f_{1}(\lambda) := \mu$ with $$\mu_{e} = ((\lambda_{o_{1}}+1),\dots,(\lambda_{o_{m}}+1))$$ and $$\mu_{o} = ((\lambda_{e_{1}}+1),\dots,(\lambda_{e_{k-r}}+1),(\lambda_{e_{k-r+1}}-1),\dots,(\lambda_{e_{k}}-1)).$$ Here we note that $\mu \in P(n)$ and $f_{1}$ reverses the parity of parts; i.e., for $\lambda$ with $k$ even and $m$ odd parts, we get $f_{1}(\lambda) = \mu$ with $k$ odd and $m$ even parts and $\mu \in P_{o}(n)$. Suppose for $\lambda^{'} \neq \lambda^{''} (\in G_{e,0}(n))$ with $\ell(\lambda^{'}) = \ell(\lambda^{''})$, we have $\mu^{'} = f_{1}(\lambda^{'}) = f_{1}(\lambda^{''}) = \mu^{''}$. Then $\ell_{e}(\lambda) = \ell_{e}(\lambda^{''})$ and so, $\ell_{o}(\lambda) = \ell_{o}(\lambda^{''})$. Now, since $\lambda^{'}$ and $\lambda^{''}$ being distinct, by the definition of $f_{1}$ we have at least a tuple $(i,j) \in \mathbb{Z}_{>0} \times \mathbb{Z}_{>0}$ such that $\mu^{'}_{i} \neq \mu^{''}_{j}$. Next, we consider the case when $\lambda^{'} \neq \lambda^{''} (\in G_{e,0}(n))$ with $\ell(\lambda^{'}) \neq \ell(\lambda^{''})$ and it is immediate that $\ell(\mu^{'}) \neq \ell(\mu^{''})$ and therefore, $\mu^{'} \neq \mu^{''}$. So, $f_{1}$ is an injective map.

For $\lambda \in G_{e,1}(n)$ with $\lambda_{e} = (\lambda_{e_{1}},\dots,\lambda_{e_{k}})$  and $\lambda_{o} = (\lambda_{o_{1}},\dots,\lambda_{o_{m}})$ where $k+m = \ell(\lambda) = 2r+1$ for some $r \in \mathbb{Z}_{> 0}$. Here we note that, $k > r$ and $k-m \geq 3$ but $k-m =1$ holds only when $a(\lambda)$ is odd, because $k = m+1$ and $a(\lambda)$ is even implies that $\lambda \in G^{0}_{e}(n)$. We exclude the condition $k = m+2$ as it contradicts that $k+m = 2r+1$.

 We define $f_{2} : G_{e,1}(n) \rightarrow P_{o}(n)$ with $f_{2}(\lambda) := \mu$, where
 \begin{equation}\label{map1}
 \mu_{e} =
  \begin{cases}
 ((\lambda_{o_{1}}+1),\dots,(\lambda_{o_{m}}+1)) \cup (\lambda_{e_{1}}+2) &\quad \text{if $a(\lambda)$ is even},\\
 ((\lambda_{o_{2}}+1),\dots,(\lambda_{o_{m}}+1)) &\quad \text{if $a(\lambda)$ is odd},\\
 \end{cases}
 \end{equation}
and
  \begin{equation}\label{map2}
 \mu_{o} =
 \begin{cases}
((\lambda_{e_{2}}+1),\dots,(\lambda_{e_{k-r-1}}+1),(\lambda_{e_{k-r}}-1),\dots,(\lambda_{e_{k}}-1))  &\quad \text{if $a(\lambda)$ is even},\\
((\lambda_{e_{1}}+1),\dots,(\lambda_{e_{k-r-1}}+1),(\lambda_{e_{k-r}}-1),\dots,(\lambda_{e_{k}}-1)) \cup (\lambda_{o_{1}}+2) &\quad \text{if $a(\lambda)$ is odd}.\\
 \end{cases}
 \end{equation}
For $a(\lambda)$ even, $$\ell_{o}(\mu)-\ell_{e}(\mu) = k-1-(m+1) =k-m-2 \geq 1$$ and for $a(\lambda)$ odd, $$\ell_{o}(\mu)-\ell_{e}(\mu) = k+1-(m-1) =k-m+2 \geq 3.$$ Hence, $\mu \in P_{o}(n)$ and by similar argument as give before, one can show that $f_{2}$ is injective.

Next, for $\lambda \in \overline{G^{0}_{e}}(n)$ with $\lambda = (\lambda_{1},\dots,\lambda_{\ell})$, we define $f_{3} : \overline{G^{0}_{e}}(n) \rightarrow P_{o}(n)$ by  $$f_{3}(\lambda) = \mu = ((\lambda_{1}+1),\lambda_{4},\dots,\lambda_{\ell}) \cup ((\lambda_{2}-2),(\lambda_{3}-2)) \cup (2,1).$$ Independent of whether $\lambda_{2}$ and $\lambda_{3}$ are odd or even, we can observe that $\ell_{o}(\mu)-\ell_{e}(\mu) = 1$ and $a(\mu) = \lambda_{1}+1$ is odd. By definition, $f_{3}$ is an injective map. Next, we show that images of $\{f_{i}\}_{1 \leq i \leq 3}$ are mutually disjoint by considering the following cases
\begin{enumerate}
	\item  By definition of the maps given before, $f_{1}(G_{e,0}(n)) \varsubsetneq P^{0}_{o}(n)$ where $$P^{0}_{o}(n) := \{\mu \in P_{o}(n) : \ell(\mu) \ \equiv 0 \ ( \text{mod} \ 2) \}$$ and $f_{2}(G_{e,1}(n)) \varsubsetneq P^{1}_{o}(n)$ with $$P^{1}_{o}(n) := \{\mu \in P_{o}(n) : \ell(\mu) \ \equiv 1 \ ( \text{mod} \ 2) \}.$$ So, $f_{1}(G_{e,0}(n)) \ \cap \ f_{2}(G_{e,1}(n)) = \emptyset$.
	\item For $\lambda \in G_{e,0}(n)$ with $\ell_{e}(\lambda)-\ell_{o}(\lambda) \geq 2$, we have $f_{1}(\lambda) = \mu \in P_{o}(n)$ with $\ell_{o}(\mu)-\ell_{e}(\mu) \geq 2$ and for $\lambda \in G_{e,0}(n)$ with $\ell_{e}(\lambda)-\ell_{o}(\lambda) = 1$, $f_{1}(\lambda) = \mu \in P_{o}(n)$ with $\ell_{o}(\mu)-\ell_{e}(\mu) = 1$ but then $a(\mu) = \lambda_{1}+1$ is even. Considering $\lambda \in \overline{G^{0}_{e}}(n)$, $f_{3}(\lambda) = \mu \in P_{o}(n)$ with $a(\mu) = \mu_{1}$ is odd and $\ell_{o}(\mu)-\ell_{e}(\mu) = 1$. Therefore, $f_{1}( G_{e,0}(n)) \ \cap \ f_{3}(\overline{G^{0}_{e}}(n)) = \emptyset$.
	\item Let us consider $\lambda \in G_{e,1}(n)$ with $a(\lambda)$ is even and $\ell_{e}(\lambda)-\ell_{o}(\lambda)=3$. Then $f_{2}(\lambda) = \mu \in P_{o}(n)$ with $a(\mu)$ is even and $\ell_{o}(\mu)-\ell_{e}(\mu)=1$. For $\ell_{e}(\lambda)-\ell_{o}(\lambda) \geq 4$, we have $f_{2}(\lambda) = \mu$ with $\ell_{o}(\mu)-\ell_{e}(\mu) \geq 2$. Moreover, if $a(\lambda)$ is odd then it is immediate that $\ell_{o}(\mu)-\ell_{e}(\mu) \geq 3$ and consequently, $f_{2}( G_{e,1}(n)) \ \cap \ f_{3}(\overline{G^{0}_{e}}(n)) = \emptyset$. So, the map $f : \overline{G_{e}}(n) \rightarrow P_{o}(n)$ is injective.
\end{enumerate}
For $\mu \in P_{o}(n)$ with its odd component $\mu_{o} =(\mu_{o_{1}},\dots,\mu_{o_{s}})$, we define
\begin{equation*}
\overline{P_{o}}(n) := \{\mu \in P_{o}(n) : \ell_{e}(\mu) =2 \ \text{and} \ \mu_{o_{i}}=1 \ \text{for all} \ 1 \leq i \leq s \}.
\end{equation*}
By the definition of $f$, it is clear that $\overline{P_{o}}(n) \varsubsetneq P_{o}(n) \setminus f(\overline{G_{e}}(n))$. Now, it remains to show that $|\overline{P_{o}}(n)| > |P_{e}(n) \setminus \overline{G_{e}}(n)|$.

For $n$ even and for $\lambda \in \overline{P_{o}}(n)$, we have $\ell_{o}(\lambda) = 2k+2$ for some $k \in \mathbb{Z}_{>0}$. Here we observe that
\begin{equation}\label{card1}
|\{\lambda \in P(n) : \lambda_{1}+\lambda_{2}=n; \lambda_{1}, \lambda_{2} \ \text{both even}                     \}|=\Bigl\lfloor\frac{n}{4}\Bigr\rfloor
\end{equation}
and 
\begin{equation}\label{card2}
|\{\lambda \in P(n) : \lambda_{1}+\lambda_{2}=n; \lambda_{1} \ \text{even}, \lambda_{2} \ \text{odd and} \ \lambda_{2} \in \mathbb{Z}_{\geq 3}                     \}|=\Bigl\lfloor\frac{n-3}{4}\Bigr\rfloor.
\end{equation}
Since, $\lambda \in \overline{P_{o}}(n)$ with $n$ even positive integer and $\ell_{o}(\lambda) = 2k+2$, for each $k \in \mathbb{Z}_{>0}$, then by \eqref{card1}, 
\begin{equation}\label{eq1}
|\{\lambda \in \overline{P_{o}}(n) : \lambda_{1}+\lambda_{2}+(2k+2) \times 1=n; \lambda_{1}, \lambda_{2} \ \text{both even}                     \}|=\Bigl\lfloor\frac{n-2k-2}{4}\Bigr\rfloor
\end{equation}
and $1\leq k \leq \dfrac{n-6}{2}$ because $k$ maximizes only when both $\lambda_{1}$ and $\lambda_{2}$ minimum; i.e., only the instance $2+2+(2k+2)\times1=n$ which implies $k = \dfrac{n-6}{2}$. Therefore we have,
\begin{equation}\label{eq2}
|\overline{P_{o}}(n)| = \sum_{k=1}^{\frac{n-6}{2}}\Bigl\lfloor\frac{n-2k-2}{4}\Bigr\rfloor.
\end{equation}

Similarly, for $n$ odd, we have $\ell_{o}(\lambda) = 2k+1$ with $1 \leq k \leq \dfrac{n-5}{2}$ and 
\begin{equation}\label{eq3}
|\{\lambda \in \overline{P_{o}}(n) : \lambda_{1}+\lambda_{2}+(2k+1) \times 1=n; \lambda_{1}, \lambda_{2} \ \text{even}                     \}|=\Bigl\lfloor\frac{n-2k-1}{4}\Bigr\rfloor.
\end{equation}
Consequently, 
\begin{equation}\label{eq4}
|\overline{P_{o}}(n)| = \sum_{k=1}^{\frac{n-5}{2}}\Bigl\lfloor\frac{n-2k-1}{4}\Bigr\rfloor.
\end{equation}

Now for $n$ even we will show that
\begin{equation}\label{eq5}
|P_{e}(n) \setminus \overline{G_{e}}(n)| = 1 + \sum_{k=1}^{\lfloor \frac{n-2}{6}\rfloor}\Bigl\lfloor \frac{n-6k+2}{4}\Bigr\rfloor + \sum_{k=1}^{\lfloor \frac{n-6}{6}\rfloor}\Bigl\lfloor \frac{n-6k-2}{4}\Bigr\rfloor.
\end{equation}
We interpret the set $P_{e}(n) \setminus \overline{G_{e}}(n)$ as a disjoint union of its three proper subsets given by\\
$P_{e}(n) \setminus \overline{G_{e}}(n) = A_{1} \cup A_{2} \cup A_{3}$ where,
\begin{equation}
\begin{split}
A_{1} & = \{\lambda \in P_{e}(n) \setminus \overline{G_{e}}(n): 0 \leq \lambda_{3} \leq 1\},\\
A_{2} & = \bigcup_{k \geq 1} A_{2,k}, \\
A_{3} & = \bigcup_{k \geq 1} A_{3,k}, \\
&\text{with}\\
A_{2,k} & = \{\lambda=(\lambda_1,\lambda_2,2,\dots,2,1,\dots,1)\vdash n: \lambda_1\ \text{and}\ \lambda_2\ \text{even},   \mult_{\lambda}(2)=2k-1, \mult_{\lambda}(1)=2k \},\\
A_{3,k} & = \{\lambda=(\lambda_1,\lambda_2,2,\dots,2,1,\dots,1)\vdash n: \lambda_1\ \text{even},\lambda_2\ \text{odd},   \mult_{\lambda}(2)=2k, \mult_{\lambda}(1)=2k-1\}.
\end{split}
\end{equation}

Next, we explicitly describe the sets and will derive their cardinality by separating into three cases.

\textit{Case 1(E):} We observe that $|A_{1}|=1$ because we have only one possibility $(\lambda_{1},\lambda_{2},\lambda_{3})=(\lambda_{1},0,0)$. We reject the other three possibilities; i.e., $(\lambda_{1},\lambda_{2},\lambda_{3})=(\lambda_{1},0,1)$ as $\lambda_{2} \geq \lambda_{3}$, $(\lambda_{1},\lambda_{2},\lambda_{3})=(\lambda_{1},1,0)$ as $n$ even and $(\lambda_{1},\lambda_{2},\lambda_{3})=(\lambda_{1},1,1)$ as $\lambda \in P_{e}(n) \setminus \overline{G_{e}}(n)$. Next, we look at the subset of $A_1$, say $A_{1,\geq 2}:=\{\lambda\in A_1: \lambda_2 \geq 2\}$ and note that $A_{1,\geq 2}=0\emptyset$. This is because for $\lambda \in A_{1,\geq 2}$, there are altogether four possibilities for $\lambda_3 \in \{0,1\}$.\\ 
For $\lambda_3=0$, the choice $(\lambda_1,\lambda_2,\lambda_3)=(\lambda_{1},\lambda_{2},0)$ and $\lambda_2$ is even is impossible as $\lambda \in P_{e}(n) \setminus \overline{G_{e}}(n)$ and if $\lambda_2$ is odd, again an impossible option since $n$ is even. Whereas for $\lambda_3=1$, the choice $(\lambda_1,\lambda_2,\lambda_3)=(\lambda_{1},\lambda_{2},1)$ and $\lambda_2$ is even is impossible as $n$ is even and if $\lambda_2$ is odd, again an impossible option since $\lambda \in P_{e}(n) \setminus \overline{G_{e}}(n)$.

\textit{Case 2(E):} By \eqref{card1},
\begin{equation}\label{eq6}
|A_{2,k}|=\Bigl \lfloor\frac{n-6k+2}{4}\Bigr\rfloor
\end{equation}
and $1\leq k\leq\lfloor\frac{n-2}{6}\rfloor$ because $k$ maximizes only when both $\lambda_{1} \ \text{and} \ \lambda_{2}$ 
minimum; i.e., the instance $2+2+(2k-1)\times 2+(2k)\times 1=n$ which implies $k \leq \lfloor\frac{n-2}{6}\rfloor$. By \eqref{eq6},
\begin{equation}\label{req1}
A_{2}=\bigcup_{k=1}^{\lfloor\frac{n-2}{6}\rfloor} A_{2,k} \ \text{and} \ |A_{2}|=\sum_{k=1}^{\lfloor\frac{n-2}{6}\rfloor} \Bigl \lfloor\frac{n-6k+2}{4}\Bigr\rfloor.
\end{equation} 

\textit{Case 3(E):} From \eqref{card2}, it follows that
\begin{equation}\label{eq7}
|A_{3,k}|=\Bigl \lfloor\frac{(n-6k+1)-3}{4}\Bigr\rfloor=\Bigl \lfloor\frac{n-6k-2}{4}\Bigr\rfloor
\end{equation}
and $1\leq k\leq\lfloor\frac{n-6}{6}\rfloor$ because $k$ maximizes only when both $\lambda_{1} \ \text{and} \ \lambda_{2}$ 
minimum; i.e., the instance $4+3+(2k)\times 2+(2k-1)\times 1=n$ which implies $k \leq \lfloor\frac{n-6}{6}\rfloor$. By \eqref{eq7},
\begin{equation}\label{req2}
A_{3}=\bigcup_{k=1}^{\lfloor\frac{n-6}{6}\rfloor} A_{3,k} \ \text{and} \ |A_{3}|=\sum_{k=1}^{\lfloor\frac{n-6}{6}\rfloor} \Bigl \lfloor\frac{n-6k-2}{4}\Bigr\rfloor.
\end{equation}
By \textit{Case 1(E)}, \eqref{req1} and \eqref{req2} we have \eqref{eq5}.

For all $n$ odd integers greater equal $9$, we will show that 
\begin{equation}\label{eq8}
|P_{e}(n) \setminus \overline{G_{e}}(n)| = 1+ \Bigl\lfloor \frac{n-4}{4}\Bigr\rfloor+ \sum_{k=1}^{\lfloor \frac{n-5}{6}\rfloor}\Bigl\lfloor \frac{n-6k-1}{4}\Bigr\rfloor + \sum_{k=1}^{\lfloor \frac{n-9}{6}\rfloor}\Bigl\lfloor \frac{n-6k-5}{4}\Bigr\rfloor.
\end{equation}
Similarly as before, we write $P_{e}(n) \setminus \overline{G_{e}}(n)$ as a disjoint union of its four proper subsets given by
$P_{e}(n) \setminus \overline{G_{e}}(n) =B_{0} \cup B_{1} \cup B_{2} \cup B_{3}$ where,
\begin{equation}
\begin{split}
B_{0} &= \{\lambda=(\lambda_1,\lambda_2,1) \in P_{e}(n) \setminus \overline{G_{e}}(n):\lambda_{2} \geq 4 \}\\
B_{1} & = \{\lambda \in P_{e}(n) \setminus \overline{G_{e}}(n): 0 \leq \lambda_{2} \leq 2 \ \text{and} \ 0 \leq \lambda_{3} \leq 1 \}\\
B_{2} & = \bigcup_{k \geq 1} B_{2,k} \\
B_{3} & = \bigcup_{k \geq 1} B_{3,k}, \\
&\text{with}\\
B_{2,k} & = \{\lambda=(\lambda_1,\lambda_2,2,\dots,2,1,\dots,1)\vdash n: \lambda_1\ \text{and}\ \lambda_2\ \text{even},   \mult_{\lambda}(2)=2k, \mult_{\lambda}(1)=2k+1 \},\\
B_{3,k} & = \{\lambda=(\lambda_1,\lambda_2,2,\dots,2,1,\dots,1)\vdash n: \lambda_1\ \text{even},\lambda_2\ \text{odd},   \mult_{\lambda}(2)=2k+1, \mult_{\lambda}(1)=2k\}.
\end{split}
\end{equation}
\textit{Case 1(O):} For $\lambda= (\lambda_1,\lambda_2,1) \in B_0$ and $n$ is odd, it follows that both $\lambda_1$ and $\lambda_2$ are even. Therefore minimal choice for $n$ is $9$ because otherwise $\lambda_1 \geq \lambda_2 \geq 4$ with the constraint that both $\lambda_1$ and $\lambda_{2}$ even would be an impossibility in such context. Moreover, we can observe that 
\begin{equation}\label{new0}
|B_0| =
\Bigl\lfloor \frac{n-4}{4}\Bigr\rfloor.
\end{equation}

\textit{Case 2(O):} We observe that $|B_{1}|=1$ because we have only one possibility $(\lambda_{1},\lambda_{2},\lambda_{3})=(\lambda_{1},2,1)$. We reject the other three possibilities; i.e., $(\lambda_{1},\lambda_{2},\lambda_{3})=(\lambda_{1},0,1)$ as $\lambda_{2} \geq \lambda_{3}$, $(\lambda_{1},\lambda_{2},\lambda_{3})=(\lambda_{1},0,0)$ and $(\lambda_{1},\lambda_{2},\lambda_{3})=(\lambda_{1},2,0)$ as $n$ odd, $(\lambda_{1},\lambda_{2},\lambda_{3})=(\lambda_{1},1,0)$ and $(\lambda_{1},\lambda_{2},\lambda_{3})=(\lambda_{1},1,1)$ as $\lambda \in P_{e}(n) \setminus \overline{G_{e}}(n)$.

\textit{Case 3(O):} By \eqref{card1},
\begin{equation}\label{eq9}
|B_{2,k}|=\Bigl \lfloor\frac{n-6k-1}{4}\Bigr\rfloor
\end{equation}
and $1\leq k\leq\lfloor\frac{n-5}{6}\rfloor$ because $k$ maximizes only when both $\lambda_{1} \ \text{and} \ \lambda_{2}$ 
minimum; i.e., the instance $2+2+(2k)\times 2+(2k+1)\times 1=n$ which implies $k \leq \lfloor\frac{n-5}{6}\rfloor$. By \eqref{eq9},
\begin{equation}\label{req3}
B_{2}=\bigcup_{k=1}^{\lfloor\frac{n-5}{6}\rfloor} B_{2,k} \ \text{and} \ |B_{2}|=\sum_{k=1}^{\lfloor\frac{n-5}{6}\rfloor} \Bigl \lfloor\frac{n-6k-1}{4}\Bigr\rfloor.
\end{equation}

\textit{Case 4(O):} From \eqref{card2}, it follows that
\begin{equation}\label{eq10}
|B_{3,k}|=\Bigl \lfloor\frac{(n-6k-2)-3}{4}\Bigr\rfloor=\Bigl \lfloor\frac{n-6k-5}{4}\Bigr\rfloor
\end{equation}
and $1\leq k\leq\lfloor\frac{n-9}{6}\rfloor$ because $k$ maximizes only when both $\lambda_{1} \ \text{and} \ \lambda_{2}$ 
minimum; i.e., the instance $4+3+(2k+1)\times 2+(2k)\times 1=n$ which implies $k \leq \lfloor\frac{n-9}{6}\rfloor$. By \eqref{eq10},
\begin{equation}\label{req4}
B_{3}=\bigcup_{k=1}^{\lfloor\frac{n-9}{6}\rfloor} B_{3,k} \ \text{and} \ |B_{3}|=\sum_{k=1}^{\lfloor\frac{n-9}{6}\rfloor} \Bigl \lfloor\frac{n-6k-5}{4}\Bigr\rfloor.
\end{equation}
By \textit{Case 2(O)}, \eqref{new0}, \eqref{req3} and \eqref{req4} we have \eqref{eq8}.

Therefore, by Lemmas \ref{lem1} and \ref{lem2}, $|\overline{P_{o}}(n)| > |P_{e}(n) \setminus \overline{G_{e}}(n)|$ for all $n \geq \mathbb{Z}_{\geq 14} \cup \{11,13\}$. To conclude the proof, it remains to check for $n \in \{1,3,4,5,6,7,8,9,10,12\}$ which we did by numerically checking in Mathematica.

\section{Proof of Theorem \ref{thm:distinct}}\label{sec:distinct}

Following the definitions in Section \ref{sec:proof-kkl}, set
\begin{equation*}
\begin{split}
H^0_e(n) &:= G^0_e(n) \cap D(n),\\
\overline{H^{0}_{e}}(n) &:= \{\lambda \in H^0_e(n): \ell_o(\lambda)>1\},\\
\text{and}\ \ H_e(n)&:= G_e(n) \cap D(n).
\end{split}
\end{equation*}
We split $H_e(n)$ into $H_{e,0}(n)=G_{e,0}(n) \cap D(n)$ and $H_{e,1}(n)=G_{e,1}(n) \cap D(n)$. Similarly, define the map $f: H_e(n) \rightarrow D_o(n)$ by $f|_{H_{e,0}(n)}=f_1$ and $f|_{H_{e,1}(n)}=f_2$. Since $H_e(n) \subsetneq G_e(n)$, we conclude that the map $f$ is injective by \eqref{map1} and \eqref{map2}.  

Now we are to show that $d_o(n) - |f(H_e(n))|>d_e(n) - |H_e(n)|$ for all $n>31$. The subset $D_o(n) \setminus f(H_e(n))$ contains different classes of partitions. One of which is 
$$\overline{D_{o}}(n):= \{\lambda \in D_o(n) \setminus f(H_e(n)) : \ell_o(\lambda)-\ell_e(\lambda)=1\ \text{and}\ a(\lambda) \equiv 1 (\text{mod}\ 2)\}.$$ We note that $\overline{D_{o}}(n)$ may contain other classes of partitions depending on $n$ is even or odd.

For a partition $\lambda \in \overline{H^{0}_{e}}(n)$, we split $\lambda$ into its even component $\lambda_e=(\lambda_{e_{1}},\lambda_{e_{2}},\dots,\lambda_{e_{m+1}})$ and $\lambda_o=(\lambda_{o_{1}},\lambda_{o_{2}},\dots,\lambda_{o_{m}})$ for some $m \in \mathbb{Z}_{\geq 2}$. Now, me make a transformation of $\lambda$ into $\lambda^\ast$ with
\begin{equation*}
\lambda^\ast = (\lambda_{e_{1}}+\lambda_{o_{1}},\lambda_{e_{2}}+\lambda_{o_{2}}) \cup (\lambda_{e_{3}},\lambda_{e_{4}},\dots,\lambda_{e_{m+1}}) \cup (\lambda_{o_{3}},\lambda_{o_{4}},\dots,\lambda_{o_{m}}) \in \overline{D_{o}}(n).
\end{equation*} 
We observe that two partitions, say $\lambda, \overline{\lambda} \in \overline{H^{0}_{e}}(n)$, where
\begin{equation*}
\begin{split}
\lambda &= (\lambda_{e_{1}},\lambda_{e_{2}},\dots,\lambda_{e_{m+1}}) \cup (\lambda_{o_{1}},\lambda_{o_{2}},\dots,\lambda_{o_{m}})\\
\text{and}\ \ \overline{\lambda} &= (\overline{\lambda}_{e_1},\overline{\lambda}_{e_2},\dots,\overline{\lambda}_{e_{m+1}}) \cup (\overline{\lambda}_{o_1},\overline{\lambda}_{o_2},\dots,\overline{\lambda}_{o_m}),
\end{split}
\end{equation*}
transform to a same partition, say $\mu \in \overline{D_{o}}(n)$ if and only if
\begin{equation*}
\lambda_{e_{1}}-\overline{\lambda}_{e_{1}} = \overline{\lambda}_{o_{1}}-\lambda_{o_{1}} \equiv 0 (\text{mod}\ 2)\ \ \text{or}/ \text{and}\ \ \lambda_{e_{2}}-\overline{\lambda}_{e_{2}} = \overline{\lambda}_{o_{2}}-\lambda_{o_{2}} \equiv 0 (\text{mod}\ 2).
\end{equation*}
If those cases arise we subtract some multiple of $2$ from the greatest part of the resultant partition and add the multiple of $2$ to the other even parts which are present in the partition, and continue this process till we have a repetition among the parts of the partition. This process is injective by its definition, and we denote it by $g$. For example, consider partitions $\lambda=(12,10,6,2)\cup(7,3,1)$ and $\overline{\lambda}=(10,8,6,2)\cup(9,5,1)$ in $\overline{H^{0}_{e}}(41)$, then both $\lambda$ and $\overline{\lambda}$ maps to the same partition $\mu=(19,13,6,2,1) \in \overline{D_{o}}(41)$. Consequently, by the process $g$, finally $\lambda\mapsto (19,13,6,2,1)$ whereas $\overline{\lambda} \mapsto (17,13,8,2,1)$. As a trivial remark, $\overline{H^{0}_{e}}(n) = \emptyset$ for all positive even integers $n \leq 14$, since $6+4+2+3+1=16$ is the least possible option. 

Depending on the parity of $n$, it remains to analyze the left over set
\begin{equation}\label{leftoverset}
\widetilde{H^0_e}(n):= \{\lambda \in H^0_e(n) : \ell_{e}(\lambda)-\ell_{o}(\lambda)=1, \ell_{o}(\lambda)\leq 1\ \text{and}\ a(\lambda) \equiv 0 (\text{mod}\ 2)\},
\end{equation}
which is unmapped yet (after applying the map $f$ and $g$).

For $n$ to be an even positive integer, we observe that $\widetilde{H^0_e}(n)$ consists of only one partition $(n)$. An even integer $n$ can be expressed as a sum of two consecutive odd integers if and only if $n$ is divisible by 4. If $n$ is divisible by 4, then for some definite odd integer $\lambda_{o_1}$ we get $(\lambda_{o_1},\lambda_{o_1}-2)\in D_o(n)$, which is not mapped yet. So we map $(n)$ to $(\lambda_{o_1},\lambda_{o_1}-2)$. If $n$ is not divisible by 4, then for some definite odd integer $\lambda_{o_1}$ we get $(\lambda_{o_1},\lambda_{o_1}-2,2)\in D_o(n)$, which is not mapped yet. So in this case we map $(n)$ to $(\lambda_{o_1},\lambda_{o_1}-2,2)$. Therefore, by some elementary observations we get that the theorem is true for all even integer $n>6$.

Let $n$ be odd. We rewrite \eqref{leftoverset} as
\begin{equation*}
\widetilde{H^0_e}(n):= \{\lambda \in H^0_e(n) : \ell_{e}(\lambda)=2, \ell_{o}(\lambda)= 1\ \text{and}\ a(\lambda) \equiv 0 (\text{mod}\ 2)\}.
\end{equation*}
Write a partition $\lambda \in \widetilde{H^0_e}(n)$ into its even component $\lambda_e=(\lambda_{e_1},\lambda_{e_{2}})$ and odd component $\lambda_o=(\lambda_{o_1})$. We split $\widetilde{H^0_e}(n)$ into following three classes:
\begin{enumerate}
	\item $\widetilde{H}^0_{e,1}(n):= \{\lambda \in \widetilde{H^0_e}(n): \lambda_{e_2}=2\}$,
	\item $\widetilde{H}^0_{e,2}(n):= \{\lambda \in \widetilde{H^0_e}(n): \lambda_{e_2}\geq 6\}$, and
	\item $\widetilde{H}^0_{e,3}(n):= \{\lambda \in \widetilde{H^0_e}(n): \lambda_{e_2}= 4\}$.
\end{enumerate}
Now we consider the following three classes of partitions from the set of partitions, say $\widetilde{D_o}(n) \subsetneq D_o(n)$ which have no preimage yet:
\begin{enumerate}
	\item $\widetilde{D}_{o,1}(n):= \{\pi \in \widetilde{D_o}(n) : \ell(\pi)=4\ \text{and}\ \pi_{o_1}-\pi_{o_2}=2 \}$,
	\item $\widetilde{D}_{o,2}(n):= \{\pi \in \widetilde{D_o}(n) : \ell_o(\pi)=3\ \text{and}\ \pi_{o_1}-\pi_{o_2}=2 \}$, and
	\item $\widetilde{D}_{o,3}(n):= \{\pi \in \widetilde{D_o}(n) : \ell_o(\pi)-\ell_{e}(\pi)=1, a(\pi)\equiv 0 (\text{mod}\ 2)\ \text{and}\ \pi_{e_1}-\pi_{o_1}=1\ \text{or}\ 3 \}$.
\end{enumerate}
Now we construct an injective map from $\widetilde{H}^0_{e,1}(n)$ to $\widetilde{D}_{o,1}(n)$. Let $\lambda=(\lambda_{e_1},2)\cup (\lambda_{o_1}) \in \widetilde{H}^0_{e,1}(n)$. Define a transformation $S$ such that $S(\lambda)=(\lambda_{e_1})\cup(\lambda_{o_1}+1,1)$. Now define $S^\ast$ such that
\begin{equation*}
S^\ast(S(\lambda))=
\begin{cases}
S(\lambda) & \quad \text{if $\lambda_{e_1} \equiv 0 (\text{mod}\ 4)$},\\
(\lambda_{e_1}-2)\cup (\lambda_{o_1}+1,3)& \quad \text{if $\lambda_{e_1} \equiv 2 (\text{mod}\ 4)$.}
\end{cases}	
\end{equation*}
Now define $S^{\ast\ast}$ such that $S^{\ast\ast}(S^\ast(S(\lambda)))=(\lambda_{o_2},\lambda_{o_2}-2)\cup(\lambda_{o_3})\cup(\lambda_{o_1}+1)$, where $\lambda_{o_2}+\lambda_{o_2}-2=\lambda_{e_1} \text{~or~} \lambda_{e_1}-2$, and $\lambda_{o_3}=1 \text{~or~} 3$ accordingly. For example: $(24,5,2)$ maps to $(13,11,6,1)$ and $(22,7,2)$ maps to $(11,9,8,3)$. This process is injective.

Our next objective is to embed the set $\widetilde{H}^0_{e,2}(n)$ into a subset of $\widetilde{D_o}(n)$ which is not mapped till now. Define a transformation $U$ such that $U(\lambda)=((\lambda_{e_1-3},3)\cup (\lambda_{o_{1}}))\cup(\lambda_{e_{2}},2)$ for $\lambda \in \widetilde{H}^0_{e,2}(n)$. Associated with $U$, let us define $U^\ast$ in such a way that
\begin{equation*}
U^\ast(U(\lambda))=
\begin{cases}
U(\lambda) &\quad~\text{if}~ \lambda_{o_1}\neq3 ~\text{and}~ \lambda_{e_1}-3\neq \lambda_{o_1},\\
(\lambda_{e_1}-3,5,1)\cup(\lambda_{e_2}-2,2) &\quad~\text{if}~ \lambda_{o_1}=3,\\
\bigl(\text{this transformation is impossible for}\ n \leq 17\bigr)\\
\bigl(\text{e.g.}\ (10,6,3)\mapsto (7,5,4,2,1)\bigr)\\
((\lambda_{e_1}-3,\lambda_{o_1}-2)\cup (5))\cup(\lambda_{e_2}-2,2) &\quad ~\text{if}~ \lambda_{e_1}-3=\lambda_{o_1},\\
\bigl(\text{this transformation is impossible for}\ n \leq 25\bigr) \\
\bigl(\text{e.g.}\ (12,9,6)\mapsto (9,7,5,4,2)\bigr)\\
((\lambda_{e_1}-3,\lambda_{o_1}-4)\cup(5))\cup(\lambda_{e_2}-2,4) &\quad ~\text{if}~ \lambda_{e_1}-1=\lambda_{o_1}, \lambda_{e_2}\neq6,\\
\bigl(\text{this transformation is impossible for}\ n \leq 29\bigr)\\
\bigl(\text{e.g.}\ (12,11,8)\mapsto (9,7,6,5,4)\bigr)\\
((\lambda_{e_1}-3,\lambda_{o_1}-4)\cup(3))\cup(6,4) &\quad~\text{if}~ \lambda_{e_1}-1=\lambda_{o_1},~\text{and}~ \lambda_{e_2}=6.\\
\bigl(\text{this transformation is impossible for}\ n \leq 23 \bigr)\\
\bigl(\text{e.g.}\ (10,9,6)\mapsto (7,6,5,4,3)\bigr)
\end{cases}
\end{equation*}
Denote the resulting transformation $U^\ast U$ by $\widetilde{U}$. Note that $\ell(\widetilde{U}(\lambda))=5$, where the resulting partitions; i.e., images under $\widetilde{U}$, contains parts from $\{3,5\}$ and its smallest even part from $\{2,4\}$. For a partition $\lambda \in H^0_e(n)$ with $\ell(\lambda)=7$, $\ell(g(\lambda))=5$. Now, assume $\widetilde{U}(\lambda)=g(\mu)$ for some partition $\mu$ with $\ell(\mu)=7$. In the map $g$, after applying the first transformation, the other transformations are applied on the even parts only (if it necessary). If $\widetilde{U}(\lambda)=g(\mu)$, then we remove one of $2$, $3$, $4$, or $5$ (which one exists in $g(\mu)$) from $g(\mu)$ by applying similar transformation.

Now we compare the number of partitions in $\widetilde{H}^0_{e,3}(n)$ with $\widetilde{D}_{o,2}(n)$ and $\widetilde{D}_{o,3}(n)$. $$\bigl|\widetilde{H}^0_{e,3}(n)\bigr|=\left\lfloor\frac{n-3}{4}\right\rfloor.$$

Let $\pi=(\pi_{o_1},\pi_{o_1}-2,\pi_{o_3}) \in \widetilde{D}_{o,2}(n)$, and $\pi_{o_1}=2k+1$. Then the least possible value of $k$ is given by
\begin{equation*}
(2k+1)+(2k-1)\geq2\{n-(2k+1)-(2k-1)\}+6,
\end{equation*}
which implies $k=\big\lceil\frac{n+3}{6}\big\rceil$. So the total number of partitions in $\widetilde{D}_{o,2}(n)$ is at least $\big\lceil\frac{n-4k}{4}\big\rceil$, which is equal to $\big\lceil\frac{n}{4}-\big\lceil\frac{n+3}{6}\big\rceil\big\rceil$. Any odd integer $n$ is of the form $12\ell+r$, where $r=1,3,5,7,9,~\text{or}~11$. Calculating for all the six cases, we get the total number of partitions in this class is
\begin{equation*}
=\begin{cases}
\left\lfloor\frac{n}{12}\right\rfloor, ~\text{if}~ n\neq12\ell+9,\\
\left\lfloor\frac{n}{12}\right\rfloor+1, ~\text{if}~ n=12\ell+9.
\end{cases}
\end{equation*}
Now, $\left\lfloor\frac{n-3}{4}\right\rfloor-\begin{cases}
\left\lfloor\frac{n}{12}\right\rfloor, ~\text{if}~ n\neq12\ell+9,\\
\left\lfloor\frac{n}{12}\right\rfloor+1, ~\text{if}~ n=12\ell+9.
\end{cases}=\left\lfloor\frac{n}{6}\right\rfloor-1,\left\lfloor\frac{n}{6}\right\rfloor, ~\text{or}~ \left\lfloor\frac{n}{6}\right\rfloor+1.$

It remains to estimate a lower bound of $\bigl|\widetilde{D}_{o,3}(n)\bigr|$. Let $\pi=(\pi_{o_1},\pi_{o_2},\pi_{o_3})\cup(\pi_{e_1},\pi_{e_2}) \in \widetilde{D}_{o,3}(n)$. The greatest odd part $\pi_{o_1}$ is the largest possible value if $\pi$ contains $3$, $1$ and $2$. So the largest possible value of $\pi_{o_1}$ is $\frac{n-9}{2}$ or $\frac{n-7}{2}$ if $n\equiv3 \pmod 4$ or $n\equiv1 \pmod 4$ respectively. The smallest possible value of $\pi_{o_1}$ is greater than $\left\lfloor\frac{n}{5}\right\rfloor$. When $\pi_{o_1}$ is not the largest or the smallest possible value, then for each possible value of $\pi_{o_1}$, we get at least 6 partitions. So total number of partitions is greater than
\begin{align*}
& 6\times\frac{1}{2}\times\bigg\{\frac{n-9}{2}-1-\bigg(\bigg\lfloor\frac{n}{5}\bigg\rfloor+1+1\bigg)\bigg\}\\
& \geq3\times\bigg\{\frac{n-9}{2}-\frac{n}{5}-3\bigg\}\\
& =\frac{9(n-25)}{10}.
\end{align*}
Now, $\left\lfloor\frac{9(n-25)}{10}\right\rfloor > \left\lfloor\frac{n}{6}\right\rfloor+1$ for all odd positive integers $n>31$.

This proves the theorem for all $n>31$. We can verify the result numerically for $19<n\leq31$.

\section{Proof of Theorem \ref{thm:reverse}}\label{sec:reverse}
Define
\begin{equation*}
\begin{split}
g_o(\lambda)&= \ell_o(\lambda)-\ell_e(\lambda)\ \ \text{for}\ \ \lambda \in P_o(n)\\
 \text{and}\ \ g_e(\pi)&= \ell_e(\pi)-\ell_e(\pi)\ \ \text{for}\ \ \pi \in P_e(n).
\end{split}
\end{equation*}
Let us split $Q_o(n)$ into the following two disjoint subsets, defined by
\begin{equation*}
\begin{split}
I^1_o(n)&:=\{\lambda\in Q_o(n):g_o(\lambda)=1\ \text{and}\ a(\lambda) \equiv 1 (\text{mod}\ 2)\}\\
\text{and}\ \ I_o(n)&:=Q_o(n)\setminus I^1_o(n).
\end{split}
\end{equation*}
Then the map $f:I_o(n) \rightarrow Q_e(n)$ defined by the restriction maps $f_1$ and $f_2$ (cf. \eqref{map1}-\eqref{map2}) is injective and consequently, $f(I_o(n)):=\widetilde{I_e}(n) \subsetneq Q_e(n)$.

For $\lambda=(\lambda_{1},\lambda_{2},\lambda_{3},\dots,\lambda_{m}) \in I^1_o(n)$, we split $I^1_o(n)$ into 
\begin{enumerate}
	\item $I^1_{o,1}(n):=\{\lambda \in I^1_o(n):\lambda_{3}\geq 6 \}$.
	\item $I^1_{o,2}(n):=\bigcup_{t=2}^{5}I^1_{o,2_t}(n)$, where $I^1_{o,2_t}(n):=\{\lambda \in I^1_o(n):\lambda_{3}=t \}$. 
	\end{enumerate}
Define $\widetilde{I}^c_e(n):=Q_e(n)\setminus \widetilde{I_e}(n)$, and for $\pi=(\pi_1,\pi_2,\dots,\pi_s)\in \widetilde{I}^c_e(n)$ with its even component $\pi_e=(\pi_{e_1},\pi_{e_2},\dots,\pi_{e_k})$, we consider the following disjoint classes, where $n \geq 21$:
\begin{enumerate}
	\item $\widetilde{I}^c_{e,1}(n):=\{\pi \in \widetilde{I}^c_e(n): g_e(\pi)=1, a(\pi) \equiv 0 (\text{mod}\ 2), \text{and}\ \pi_1 \neq \pi_2\}$.
	\item $\widetilde{I}^c_{e,2}(n):=\{\pi \in \widetilde{I}^c_e(n): g_e(\pi)=1, a(\pi) \equiv 0 (\text{mod}\ 2), \text{and}\ \pi_1=\pi_2\}$.\\
	Example: $(4,4,3,3,3,2,2)$.
	\item $\widetilde{I}^c_{e,3}(n):=\{\pi \in \widetilde{I}^c_e(n): g_e(\pi)=1, a(\pi) \equiv 1 (\text{mod}\ 2), \text{and}\ \pi_1-\pi_2\in \{0,1\} \}$.\\ Example: $(5,5,3,2,2,2,2)$ and $(5,4,3,3,2,2,2)$.
	\item $\widetilde{I}^c_{e,4}(n):=\{\pi \in \widetilde{I}^c_e(n): \ell(\pi)=2r\ \text{and}\ \mult_{\pi}(2)\geq r+1\}$. Example: $(4,4,3,2,2,2,2,2)$.
	\item $\widetilde{I}^c_{e,5}(n):=\{\pi \in \widetilde{I}^c_e(n): \ell(\pi)=2r+1\ \text{and}\ \mult_{\pi}(2)\geq r+2\}$. Example: $(7,4,2,2,2,2,2)$.
	\item $\widetilde{I}^c_{e,6}(n):=\{\pi \in \widetilde{I}^c_e(n): g_e(\pi)=2, \ell(\pi) \equiv 0 (\text{mod}\ 2), \text{and}\ \pi_{e_1} \neq \pi_{e_2}\}$.\\ Example: $(13,8,8,7,4,2)$.
	\item $\widetilde{I}^c_{e,7}(n):=\{\pi \in \widetilde{I}^c_e(n): g_e(\pi)=3, \ell(\pi) \equiv 1 (\text{mod}\ 2), \text{and}\ \pi_{e_1}=\pi_{e_2}\}$.\\
	 Example: $(13,8,8,7,4,2,2)$.
\end{enumerate}
Now we define a map $\psi_1: I^1_{o,1}(n) \rightarrow \widetilde{I}^c_{e,1}(n)$. For $\lambda=(\lambda_{1},\lambda_{2},\lambda_{3},\dots,\lambda_{m})\in I^1_{o,1}(n)$,
\begin{equation*}
\psi_1(\lambda):=(\lambda_{1}+1,\lambda_{2}-4,\lambda_{3}-4,\lambda_{4},\dots,\lambda_{m})\cup(4,3),
\end{equation*}
and hence $\psi_1$ is a well defined injective map.

For the rest of the proof, in a partition $\lambda \vdash n$, as a part $x$ (resp. $y$) denotes $3$ or $5$ (resp. $2$ or $4$).

Now we construct injective maps on $I^1_{o,2}(n)$ by considering  the following five cases.

\textbf{Case 1:} $\lambda \in I^1_{o,2_2}(n)$.\\
This case is satisfied only when $n$ is even. For $\lambda=(\lambda_{1},\lambda_{2},2)\in  I^1_{o,2_2}(n)$, define
\begin{equation*}
\psi_2(\lambda):=\pi=(\lambda_{1}-\lambda_{2}+2,2,2,\dots,2),
\end{equation*}
with $\mult_{\pi}(2)=\lambda_{2}$. For example, $(23,5,2) \mapsto (20,2,2,2,2,2)$. One can observe that the map $\psi_2$ is injective and
\begin{equation*}
\psi_2(I^1_{o,2_2}(n)):=\widetilde{I}^1_{e,2_2}(n) \subsetneq \widetilde{I}^c_{e,4}(n),
\end{equation*}
where
\begin{equation*}
\widetilde{I}^1_{e,2_2}(n)=
\begin{cases}
\{\pi \in  \widetilde{I}^c_{e,4}(n): \pi_1 \geq 2\ \text{and}\ \pi_j=2\ \text{for all}\ j>1 \} &\quad \text{if}\ n\equiv 0 (\text{mod}\ 4),\\
\{\pi \in  \widetilde{I}^c_{e,4}(n): \pi_1 \geq 4\ \text{and}\ \pi_j=2\ \text{for all}\ j>1 \} &\quad \text{if}\ n\equiv 2 (\text{mod}\ 4).
\end{cases}
\end{equation*}

\textbf{Case 2:} $\lambda \in I^1_{o,2_3}(n)$.\\
Subdivide $I^1_{o,2_3}(n):=\bigcup_{t=1}^{3}I^1_{o,2_{3,t}}(n)$ with
\begin{equation*}
\begin{split}
I^1_{o,2_{3,1}}(n)&:=\{\lambda \in I^1_{o,2_3}(n): \lambda_2 \equiv 1(\text{mod 2})\},\\
I^1_{o,2_{3,2}}(n)&:=\{\lambda \in I^1_{o,2_3}(n): \lambda_2 \equiv 0(\text{mod 2})\ \text{and}\ \lambda_{1}\neq 5 \},\\
\text{and}\ \ \ I^1_{o,2_{3,3}}(n)&:=\{\lambda \in I^1_{o,2_3}(n): \lambda_2 \equiv 0(\text{mod 2})\ \text{and}\ \lambda_{1}= 5 \}.
\end{split}
\end{equation*}
We define the map $\psi_3:I^1_{o,2_3}(n) \rightarrow \widetilde{I}^c_{e}(n)$ by $\psi_3|_{I^1_{o,2_{3,t}}(n)}:=\psi_{3,t}$ for $1\leq t\leq 3$.

We take $\psi_{3,1}:=\psi_2$, given above and so,
\begin{multline*}
\widetilde{I}^1_{e,2_{3,1}}(n):=\psi_{3,1}(I^1_{o,2_{3,1}}(n))
\\=\{\pi \in  \widetilde{I}^c_{e,4}(n): a(\pi)=\pi_{e_1}\geq 2\ \text{or}\ 4, \mult_{\pi}(3)>1, \text{and}\ \pi_{e_j}=2\ \text{for}\ j>1\}.
\end{multline*}
For example, $(23,3,3,2,2)\mapsto (22,3,2,2,2,2)$.

We note that for any $\lambda \in I^1_{o,2_{3,2}}(n)$, $\lambda=(\lambda_1,\lambda_2,3,3,\dots,3,2,2\dots,2)$ with $\mult_{\lambda}(2)=r$ and $\mult_{\lambda}(3)=r+1$ for some $r \in \mathbb{Z}_{\geq 0}$. Define $\psi_{3,2}(\lambda):=\pi=(\lambda_{1}-4,3,3,\dots,3,2,2,\dots,2)$ with $\mult_{\pi}(3)\geq r+1$ and $\mult_{\pi}(2)=\frac{1}{2}\lambda_{2}+2$. Consequently,
\begin{multline*}
\widetilde{I}^1_{e,2_{3,2}}(n):=\psi_{3,2}(I^1_{o,2_{3,2}}(n))
\\=\{\pi \in \cup_{k=4}^5 \widetilde{I}^c_{e,k}(n): a(\pi)=\pi_{o_1}\geq 3, \mult_{\pi}(3)>1, \text{and}\ \pi_{e_i}=2\ \text{for}\ i \geq 1\}.
\end{multline*}
For example, $(21,8,3) \mapsto (17,3,2,2,2,2,2,2)$.

Similar as before, observe that for any $\lambda \in I^1_{o,2_{3,3}}(n)$, $\lambda=(5,4,3,3,\dots,3,2,2\dots,2)$ with $\mult_{\lambda}(2)=r$ and $\mult_{\lambda}(3)=r+1$ for some $r \in \mathbb{Z}_{\geq 0}$. Define $\psi_{3,3}(\lambda):=\pi=(3,3,\dots,3,2,2,\dots,2)$ with $\mult_{\pi}(3)= r+2$ and $\mult_{\pi}(2)=r+3$. Consequently,
\begin{equation*}
\widetilde{I}^1_{e,2_{3,3}}(n):=\psi_{3,3}(I^1_{o,2_{3,3}}(n))=\{\pi \in  \widetilde{I}^c_{e,3}(n): \mult_{\pi}(3)=\mult_{\lambda}(3)+1\ \text{and}\ \mult_{\pi}(2)=\mult_{\lambda}(2)+3 \}.
\end{equation*}
For example, $(5,4,3,3,2) \mapsto (3,3,3,2,2,2,2)$.

Moreover, $\bigcap_{1\leq i\leq 3}\widetilde{I}^1_{e,2_{3,i}}(n)=\emptyset$, for each $i$, $\{\widetilde{I}^1_{e,2_{3,i}}(n)\}\cap \widetilde{I}^1_{e,2_2}(n)=\emptyset$; and also $\{\widetilde{I}^1_{e,2_{3,i}}(n) \}_{1\leq i\leq 3}$ are mutually disjoint with the images under the maps $f$ and $\psi_1$ considered before.

\textbf{Case 3:} $\lambda \in I^1_{o,2_4}(n)$.\\
For each pair of $(\lambda_{1},\lambda_{2})$, which exist, we get a unique partition. This is also satisfied in the above two cases, but we have not considered this in the above two cases to find different images. Now analyze the disjoint subsets of $I^1_{o,2_4}(n)$, depending on $\lambda_2$ is even or odd, defined as follows
\begin{equation*}
\begin{split}
I^1_{o,2_{4,1}}(n)&:=\{\lambda \in I^1_{o,2_4}(n):\lambda_{2}\equiv 1 (\text{mod}\ 2)\ \text{and}\ \ell(\lambda)\leq 5\},\\
I^1_{o,2_{4,2}}(n)&:=\{\lambda \in I^1_{o,2_4}(n):\lambda_{2}\equiv 1 (\text{mod}\ 2)\ \text{and}\ \ell(\lambda)> 5\},\\
\text{and}\ \ I^1_{o,2_{4,3}}(n)&:=\{\lambda \in I^1_{o,2_4}(n):\lambda_{2}\equiv 0 (\text{mod}\ 2)\}.
\end{split}
\end{equation*}
When $\lambda \in I^1_{o,2_{4,1}}(n)$, then the partitions are of the two forms: $\lambda=(\lambda_1,\lambda_2,4)$ or $\lambda=(\lambda_{1},\lambda_{2},3)\cup(4,y)$, where $\lambda_1\geq \lambda_2\geq5$ are odds.\\
a) $\lambda_{1}>2\lambda_{2}+3$:
\begin{equation*}
(\lambda_{1},\lambda_{2},3)\cup(4,y) \xrightarrow{\psi_{4,1}}
\begin{cases}
(6,6,4)\cup\Biggl(\dfrac{\lambda_{1}+\lambda_{2}-12}{2},\dfrac{\lambda_{1}+\lambda_{2}-12}{2},3\Biggr) &\quad \text{if}\ \lambda_{1}-\lambda_{2}\equiv 0 (\text{mod}\ 4),\\
(6,6,6,y)\cup\Biggl(\dfrac{\lambda_{1}+\lambda_{2}-14}{2},\dfrac{\lambda_{1}+\lambda_{2}-14}{2},3\Biggr) &\quad \text{if}\ \lambda_{1}-\lambda_{2}\equiv 2 (\text{mod}\ 4),
\end{cases}
\end{equation*}
and
\begin{equation*}
(\lambda_{1},\lambda_{2},4) \xrightarrow{\psi_{4,1}}
\begin{cases}
(6,6,4)\cup\Biggl(\dfrac{\lambda_{1}+\lambda_{2}-12}{2},\dfrac{\lambda_{1}+\lambda_{2}-12}{2}\Biggr) &\quad \text{if}\ \lambda_{1}-\lambda_{2}\equiv 0 (\text{mod}\ 4),\\
(6,6,6)\cup\Biggl(\dfrac{\lambda_{1}+\lambda_{2}-14}{2},\dfrac{\lambda_{1}+\lambda_{2}-14}{2}\Biggr) &\quad \text{if}\ \lambda_{1}-\lambda_{2}\equiv 2 (\text{mod}\ 4).
\end{cases}
\end{equation*}
b) $2\lambda_{2}+3\geq \lambda_{1}\geq \lambda_{2}+8$:
\begin{equation*}
\begin{split}
(\lambda_{1},\lambda_{2},3)\cup(4,y) &\xrightarrow{\psi_{4,1}} (\lambda_{1}-\lambda_{2}-4,4,4)\cup(\lambda_{2},\lambda_{2},3)\\
\text{and}\ \ (\lambda_{1},\lambda_{2},4) &\xrightarrow{\psi_{4,1}} (\lambda_{1}-\lambda_{2}-4,4,4)\cup(\lambda_{2},\lambda_{2}).
\end{split}
\end{equation*}
c) $\lambda_{2}+8>\lambda_{1}$: For $\lambda_{2}>5$,
\begin{equation*}
\begin{split}
(\lambda_{1},\lambda_{2},4) &\xrightarrow{\psi_{4,1}} (\lambda_{1}-3,\lambda_{2}-3,2,2,2,2,2),\\
(\lambda_{1},\lambda_{2},4,4,3) &\xrightarrow{\psi_{4,1}} (\lambda_{1}-1,\lambda_{2}-1,3,2,2,2,2,2),\\
\text{and}\ \ (\lambda_{1},\lambda_{2},4,3,2) &\xrightarrow{\psi_{4,1}} (a-3,b-3,3,2,2,2,2,2,2).
\end{split}
\end{equation*}
Whereas for the left over partitions are mapped as follows
\begin{equation*}
\begin{split}
(7,5,4,3,2) &\xrightarrow{\psi_{4,1}} (5,4,4,4,4),\\
(9,5,4,3,2) &\xrightarrow{\psi_{4,1}} (5,4,4,4,4,2),\\
\text{and}\ \ (11,5,4,3,2) &\xrightarrow{\psi_{4,1}} (5,4,4,4,4,4).
\end{split}
\end{equation*}

For $\lambda \in I^1_{o,2_{4,2}}(n)$, can be written explicitly in the form $$\lambda=(\lambda_{1},\lambda_2,4,4,\dots,4,3,3,\dots,3,2,\dots,2),$$ with $\lambda_{1}, \lambda_{2}$ both odd, $\mult_{\lambda}(4)=s+1, \mult_{\lambda}(2)=r$, and $\mult_{\lambda}(3)=r+s$ for some $r,s \in \mathbb{Z}_{\geq 0}$ subject to the condition that $r+s\geq 2$. Then
\begin{equation*}
\lambda \xrightarrow{\psi_{4,2}}
\begin{cases}
(\lambda_{1},\lambda_{2},3,2,\dots,2) & \quad \text{if}\  r+s\equiv 1 \pmod 2 ,\\
(\lambda_{1},\lambda_{2},2,\dots,2) & \quad \text{if}\  r+s\equiv 0 \pmod 2.\\
\end{cases}
\end{equation*}
with 
\begin{equation*}
\mult_{\psi_{4,2}(\lambda)}(2)=
\begin{cases}
2s+2+3 \frac{r+s-1}{2}+r &\quad \text{if}\ r+s\equiv 1 \pmod 2,\\
2s+2+3 \frac{r+s}{2}+r &\quad \text{if}\ r+s\equiv 0 \pmod 2.
\end{cases}
\end{equation*}
For example $(15,7,4,3,3,3,2,2,2)\mapsto (15,7,3,2,2,2,2,2,2,2,2)$. By definition of the map $\psi_{4,2}$, these images are different from all the above images.

For $\lambda \in I^1_{o,2_{4,3}}(n)$, can be expressed as $\lambda=(\lambda_{1},\lambda_2,4,4,\dots,4,3,3,\dots,3,2,\dots,2)$ with $\lambda_{1}$ odd, $\lambda_{2}$ even, $\mult_{\lambda}(4)=s+1, \mult_{\lambda}(2)=r$, and $\mult_{\lambda}(3)=r+s+2$ for some $r,s \in \mathbb{Z}_{\geq 0}$. Then
\begin{equation*}
\lambda= \xrightarrow{\psi_{4,3}}
\begin{cases}
(\lambda_{1},\lambda_{2},3,2,\dots,2) & \quad \text{if}\  r+s\equiv 1 \pmod 2,\\
(\lambda_{1},\lambda_{2},2,\dots,2) & \quad \text{if}\  r+s\equiv 0 \pmod 2.\\
\end{cases}
\end{equation*}
with 
\begin{equation*}
\mult_{\psi_{4,3}(\lambda)}(2)=
\begin{cases}
2s+2+3 \frac{r+s+1}{2}+r &\quad \text{if}\ r+s\equiv 1 \pmod 2,\\
2s+2+3 \frac{r+s+2}{2}+r &\quad \text{if}\ r+s\equiv 0 \pmod 2.
\end{cases}
\end{equation*}
For example $(13,6,4,3,3,3,2)\mapsto (13,6,3,2,2,2,2,2,2)$. Moreover, the images under the map $\psi_{4,3}$ are different from all the above images.

\textbf{Case 4:} $\lambda \in I^1_{o,2_5}(n)$.\\
$\lambda=\lambda_e\cup\lambda_o \in D_e(n)$ with $\ell(\lambda)=2r=k+m$ ($m =\ell_o(\lambda)$), then $\lambda_{e_i}\neq \lambda_{e_{j}}$ for all $1\leq i\neq j \leq  r-m+1$ if $\lambda$ has a preimage under the map $f$ (defined by $f_1$ and $f_2$). More precisely, $\lambda_{e_{r-m}}\neq \lambda_{e_{r-m+1}}$. For example, $(15,12,12,12,8,8,3,3,2,2,2,2,2,2) \in D_e(n)$, but it has no preimage since $\lambda_{e_4}=\lambda_{e_5}=8$. So, when $\ell(\lambda)=(m+2)+m$, then $\lambda_{e_1}\neq \lambda_{e_2}$ in order to have a preimage under $f$. As a conclusion, we see that $\widetilde{I}^c_{e,6}(n)$ has no preimage yet.

Similarly as before, we further split $I^1_{o,2_5}(n)$ depending on $\lambda_{2}$ is even or odd, given by
\begin{equation*}
\begin{split}
I^1_{o,2_{5,1}}(n)&:=\{\lambda \in I^1_{o,2_5}(n):\lambda_{2}\equiv 1 \pmod 2\ \text{and}\ \ell(\lambda)> 5\},\\
I^1_{o,2_{5,2}}(n)&:=\{\lambda \in I^1_{o,2_5}(n):\lambda_{2}\equiv 1 \pmod 2\ \text{and}\ \ell(\lambda)= 5\},\\
\text{and}\ \ I^1_{o,2_{5,3}}(n)&:=\{\lambda \in I^1_{o,2_5}(n):\lambda_{2}\equiv 0 \pmod 2\}.
\end{split}
\end{equation*}
For $\lambda \in I^1_{o,2_{5,1}}(n)$, we want to define a map, say $\psi_{5,1}:I^1_{o,2_{5,1}}(n) \rightarrow \widetilde{I}^c_{e,6}(n)$. First, we explicitly write such $\lambda$ as
\begin{equation*}
\lambda=\lambda_e\cup\lambda_o:=(\lambda_1,\lambda_2,5,x,\dots,x)\cup(\lambda_{e_1},\lambda_{e_{2}},\dots,\lambda_{e_{r}}),
\end{equation*}
where $\lambda_{e_i}\in \{4,2\}$ for $1\leq i \leq r$ and $r\in \mathbb{Z}_{\geq 3}$ along with the property that $\lambda_{e_i}=\lambda_{e_j}$ for some $1\leq i\neq j \leq r$. Now define
\begin{equation*}
\lambda \xrightarrow{\psi_{5,1}}
\begin{cases}
(\lambda_1,x,\dots,x)\cup\Bigl(\hat{\lambda}_e\cup\Bigl(\lambda_{e_i}+\frac{\lambda_{2}+1}{2},\lambda_{e_j}+\frac{\lambda_{2}+1}{2},4\Bigr)\Bigr) & \quad \text{if}\ \lambda_{2} \equiv 3 \pmod 4\\
(\lambda_1,x,\dots,x)\cup\Bigl(\hat{\lambda}_e\cup\Bigl(\lambda_{e_i}+\frac{\lambda_{2}+3}{2},\lambda_{e_j}+\frac{\lambda_{2}+3}{2},2\Bigr)\Bigr) & \quad \text{if}\ \lambda_{2} \equiv 1 \pmod 4,
\end{cases}
\end{equation*}
where $\hat{\lambda}_e:= (\lambda_{e_{1}},\dots,\lambda_{e_{i-1}},\lambda_{e_{i+1}},\lambda_{e_{j-1}},\lambda_{e_{j+1}},\dots,\lambda_{e_r})$. For example $$(13,9,5,4,3,2,2)\mapsto(13,8,8,4,3,2).$$

The map $\psi_{5,1}$ may not be necessarily one to one. Then by adding and subtracting some multiples of 2 from the parts of the $\psi_{5,1}(\lambda)$ except the part $a(\psi_{5,1}(\lambda))$, we can map them to different partitions. For example, $(9,9,5,3,2,2,2)$ and $(9,5,5,4,4,3,2)$ both transform into $(9,8,8,3,2,2)$. Then for one preimage we can change the image to $(9,7,6,6,2,2)$.

When $\lambda \in I^1_{o,2_{5,2}}(n)$, then $n$ must be odd. Here partition $\lambda$ is of three form: $(\lambda_{1},\lambda_2,5,4,4)$, $(\lambda_{1},\lambda_2,5,4,2)$, and $(\lambda_{1},\lambda_2,5,2,2)$. For each of these partitions, define
\begin{equation*}
(\lambda_{1},\lambda_2,5,4,4) \xrightarrow{\psi_{5,2}}
\begin{cases}
(\lambda_{1}-5,\lambda_2-1,2,2,2,2,2)\cup(9) &\quad \text{if}\  \lambda_1\geq \lambda_{2}+4,\\
(\lambda_{1}-5,\lambda_2-5,2,2,2,2,2,2,2)\cup(9) &\quad \text{if}\  \lambda_1< \lambda_{2}+4,
\end{cases}
\end{equation*}
\begin{equation*}
(\lambda_{1},\lambda_2,5,4,2) \xrightarrow{\psi_{5,2}}
\begin{cases}
(\lambda_{1}-5,\lambda_2-1,2,2,2,2,2)\cup(7) &\quad \text{if}\  \lambda_1\geq \lambda_{2}+4,\\
(\lambda_{1}-5,\lambda_2-5,2,2,2,2,2,2,2)\cup(7) &\quad \text{if}\  \lambda_1< \lambda_{2}+4,
\end{cases}
\end{equation*}
and
\begin{equation*}
(\lambda_{1},\lambda_2,5,2,2) \xrightarrow{\psi_{5,2}}
\begin{cases}
(\lambda_{1}-5,\lambda_2-1,2,2,2,2,2)\cup(5) &\quad \text{if}\  \lambda_1\geq \lambda_{2}+4,\\
(\lambda_{1}-5,\lambda_2-5,2,2,2,2,2,2,2)\cup(5) &\quad \text{if}\  \lambda_1< \lambda_{2}+4.
\end{cases}
\end{equation*}
For example $(13,9,5,4,2)\mapsto(8,8,7,2,2,2,2,2)$ and $(11,9,5,4,4)\mapsto(9,6,4,2,2,2,2,2,2,2)$. Here, $\lambda_1-5,\lambda_2-5\geq4$ for all $n\geq31$. If a partition $\lambda$ has $3$ parts which are greater then equal to $3$, then $\mult_{\psi_{5,2}(\lambda)}(2)\geq 5$ so as to belong to the $\widetilde{I}^c_{e,4}(n)$; i.e., $\psi_{5,2}$ is an injective map to $\widetilde{I}^c_{e,4}(n)$, for all $n\geq31$.

Consider $\lambda \in I^1_{o,2_{5,3}}(n)$. Write $\lambda=(\lambda_1,5,x,\dots,x)\cup(\lambda_2,y,\dots,y)$ with $\lambda_1\geq 7$ odd, $\lambda_2\geq 6$ even, and $\text{total no. of}\ x=\text{total no. of}\ y=r$ for some non-negative integer $r$.\\
a) If $\lambda_1>2\lambda_2+3$, then 
\begin{equation*}
\lambda \xrightarrow{\psi_{5,3}}
\begin{cases}
(7,5,x,\dots,x)\cup\Bigl(\frac{\lambda_1-7}{2},\frac{\lambda_1-7}{2},\lambda_2,y,\dots,y\Bigr) &\quad \text{if}\ \lambda_1 \equiv 3 \pmod 4,\\
(5,5,x,\dots,x)\cup\Bigl(\frac{\lambda_1-5}{2},\frac{\lambda_1-5}{2},\lambda_2,y,\dots,y\Bigr) &\quad \text{if}\ \lambda_1 \equiv 1 \pmod 4.
\end{cases}
\end{equation*}
b) If $2\lambda_2+3\geq \lambda_{1}>\lambda_{2}+5$, then $\lambda \xrightarrow{\psi_{5,3}} (\lambda_1-\lambda_2-5,5,y,\dots,y)\cup(\lambda_2,\lambda_2,x,\dots,x,2)$.\\
c) If $\lambda_1=\lambda_2+5, ~\text{or}~ \lambda_2+3$, then define
\begin{equation*}
\lambda \xrightarrow{\psi_{5,3}}
\begin{cases}
(\lambda_2-1,5,x,\dots,x) \cup \Bigl((\lambda_2-2)\cup(6,y,\dots,y,2)\Bigr) &\quad \text{if}\ \lambda_{1}=\lambda_2+5,\\
(\lambda_2-1,5,x,\dots,x) \cup (\lambda_2-2,4,y,\dots,y,2) &\quad \text{if}\ \lambda_{1}=\lambda_2+3.
\end{cases}
\end{equation*}
Here we exclude the partitions $\lambda$ of the form $(13,5,x,\dots,x)\cup(8,y,\dots,y)$, since $(13,8,5)$ maps to $(7,6,6,5,2)$, and from (b) we see that $(15,6,5)$ also maps to $(7,6,6,5,2)$.\\
d) It remains to consider the two left-over classes of partitions, given by 
\begin{equation*}
\begin{split}
\lambda &= (13,5,x,\dots,x)\cup(8,y,\dots,y)\\
\text{and}\ \ \lambda &=(\lambda_1,5,x,\dots,x)\cup(\lambda_2,y\dots,y)\ \text{with}\ \lambda_1=\lambda_2+1.
\end{split}
\end{equation*} 
For these two classes of $\lambda$, define
\begin{equation*}
\begin{split}
 (13,5,x,\dots,x)\cup(8,y,\dots,y) &\xrightarrow{\psi_{5,3}} (x,\dots,x)\cup(12,12,y,\dots,y,2)\\
 \text{and}\ \ (\lambda_2+1,5,x,\dots,x)\cup(\lambda_2,y\dots,y) &\xrightarrow{\psi_{5,3}} (x,\dots,x)\cup(\lambda_2,\lambda_2,6,y,\dots,y).
\end{split}
\end{equation*}
If in two partitions of the form $(\lambda_2+1,5,x,\dots,x)\cup(\lambda_2,y,\dots,y)$ the part $\lambda_2$ is different, then the partitions are not equal. Again, when $\lambda_2$ is equal in two such partitions, then the partitions are different if and only if the combination of $x$'s and $y$'s are different. The same is true for the partitions of the form $(13,5,x,\dots,x)\cup(8,y,\dots,y)$. So the map gives us different images, which are in the subset $\widetilde{I}^c_{e,7}(n)$. By similar argument to the partitions of the $\widetilde{I}^c_{e,6}(n)$, we observe that $\widetilde{I}^c_{e,7}(n)$ have also no preimage under $f$.

\textbf{Case 5:} The partition $\lambda=(n)$.\\
If $n$ is odd, then the partition $\lambda \in Q_o(n)$, which has not mapped yet. Let $$\pi=(4,4,\ldots,4,3) ~\text{or}~ (4,4,\ldots,4,3,2).$$ Then $\pi \in Q_e(n)$, and it has no preimage yet. So we map $(n)$ to $\pi$.

Now listing the partitions of $n$ for $n\leq20$ and $n=21,23,25,27,29$, we get the inequality is true for all even numbers, and for all odd numbers greater than $7$.

\section{Proofs of Theorems \ref{thm:missingset1} and \ref{thm:missingset2}}\label{sec:missingset}

\emph{Proof of Theorem \ref{thm:missingset1}}:
Throughout this proof, $S:=\{2\}$. Define
\begin{equation*}
\begin{split}
P^S_{e,e}(n)&:= \{\lambda \in P^S_e(n): \ell(\lambda) \equiv 0 \pmod 2\}\\
\text{and}\ \ P^S_{e,o}(n)&:= \{\lambda \in P^S_e(n): \ell(\lambda) \equiv 1 \pmod 2\}.
\end{split} 
\end{equation*}

We split $P^S_e(n)$ into the following five disjoint classes.
\begin{enumerate}
	\item $P^S_{e,1}(n):= \{\lambda \in P^S_{e,e}(n): \lambda_i\neq 1\ \text{for all}\ i \}$,
	\item $P^S_{e,2}(n):= \{\lambda \in P^S_{e,o}(n): a(\lambda)\equiv 0 (\text{mod}\ 2), \ell_{e}(\lambda)-\ell_{o}(\lambda)\geq 2,\text{and}\ \lambda_i\neq 1\ \text{for all}\ i \}$,
	\item $P^S_{e,3}(n):= \{\lambda \in P^S_{e,e}(n): \lambda_j= 1\ \text{for some}\ j \}$,
	\item $P^S_{e,4}(n):= \{\lambda \in P^S_{e,o}(n): \lambda_j= 1\ \text{for some}\ j \}$, and
	\item $P^S_{e,5}(n):= \{\lambda \in P^S_{e,o}(n): a(\lambda)\equiv 0 (\text{mod}\ 2), \ell_{e}(\lambda)=\ell_{o}(\lambda)+1,\text{and}\ \lambda_i\neq 1\ \text{for all}\ i \}$.
\end{enumerate}

We apply the injective map $f$ on $P^S_{e,1}(n)$ and $P^S_{e,2}(n)$. For the sets $P^S_{e,3}(n)$ and $P^S_{e,4}(n)$, we apply $f$ in a slightly different way. For a partition $\lambda \in P^S_{e,3}(n) \cup P^S_{e,4}(n)$, we first split it as $\lambda=\hat{\lambda}\cup(1,1,\dots,1)$, where $\hat{\lambda}$ has no parts equal to 1. Then we define $\hat{f}(\lambda):=f(\hat{\lambda})\cup(1,\dots,1)$. For example, $(8,8,7,6,4,4,1,1) \mapsto (9,9,8,5,3,3,1,1)$, and $(8,8,7,6,4,1,1) \mapsto (10,8,7,5,3,1,1)$.

Now, we dissect the set $P^S_{e,5}(n)$ into two disjoint classes. For $\lambda \in P^S_{e,5}(n)$ with its even component (resp. odd component) $\lambda_e=(\lambda_{e_{1}},\lambda_{e_{2}},\dots,\lambda_{e_{r+1}})$ (resp. $\lambda_o=(\lambda_{o_{1}},\lambda_{o_{2}},\dots,\lambda_{o_r})$), define
\begin{enumerate}
	\item [(5a)] $\overline{P}^S_{e,5}(n):= \{\lambda \in  P^S_{e,5}(n) : \lambda_{e_1}\neq \lambda_{e_2}\ \text{and}\  \lambda_{e_{r+1}}\geq 6\}$, and
	\item [(5b)] $\overline{P}^{S,c}_{e,5}(n):= P^S_{e,5}(n) \setminus \overline{P}^S_{e,5}(n)$.
\end{enumerate}

Define a map, say $\phi:\overline{P}^S_{e,5}(n) \rightarrow P^S_o(n)$ by $\phi(\lambda):=\pi_o \cup \pi_e$, where
\begin{equation*}
\pi_o=\lambda_o\cup(1,1,\dots,1,1)\ \ \text{and}\ \ \pi_e=(\lambda_{e_1}-4,\lambda_{e_2}-2,\dots,\lambda_{e_{r+1}}-2).
\end{equation*}
By definition of $\phi$, for $\lambda \in \overline{P}^S_{e,5}(n)$ with $\ell_o(\lambda)=r$ it follows that
\begin{equation*}
\begin{split}
\phi(\overline{P}^S_{e,5}(n))&:=\overline{P}^S_{o,5}(n):=\{\pi \in P^S_o(n):\ell_{e}(\pi)=r+1,\ell_{o,1}(\pi)=r\ \text{and}\ \mult_{\pi}(1)=2r+4 \},
\end{split}
\end{equation*}
$\text{where}\ \ \ell_{o,1}(\lambda)\ \text{denote the number of  odd parts of}\ \lambda \vdash n\ \text{which are greater than}\ 1$.

 Therefore, $\overline{P}^S_{o,5}(n)$ is disjoint with an image of any partition belongs to $\bigcup^{4}_{i=1}P^S_{e,i}(n)$ under the map $f$ or $\hat{f}$. Since the odd parts of $\lambda$ are not altered, so the map is injective.
 
Next, we further split $\overline{P}^{S,c}_{e,5}(n)$ into four disjoint subsets.
\begin{enumerate}
	\item [(5b1)] $\overline{P}^{S,c}_{e,5_1}(n):=\{\lambda \in \overline{P}^{S,c}_{e,5}(n):\lambda_i \neq3\}$,
	\item [(5b2)] $\overline{P}^{S,c}_{e,5_2}(n):=\{\lambda \in \overline{P}^{S,c}_{e,5}(n): \ell(\lambda) \neq 3\ \text{and}\ \mult_{\lambda}(3)=1\}$,
	\item [(5b3)] $\overline{P}^{S,c}_{e,5_3}(n):=\{\lambda \in \overline{P}^{S,c}_{e,5}(n): \ell(\lambda) \neq 5\ \text{and}\ \mult_{\lambda}(3)\geq2\}$, and
	\item [(5b4)] $\overline{P}^{S,c}_{e,5_4}(n):=\bigcup_{k=1}^{2}\overline{P}^{S,c}_{e,5_{4,k}}(n),$ with
	\begin{equation*}
	\overline{P}^{S,c}_{e,5_{4,k}}(n)=\{\lambda \in \overline{P}^{S,c}_{e,5}(n):\ell(\lambda)= 2k+1\ \text{and}\ \mult_{\lambda}(3)=k\}.
	\end{equation*}
\end{enumerate}

We construct a map $\rho: \overline{P}^{S,c}_{e,5}(n) \rightarrow P^S_o(n)$ by defining
\begin{equation*}
\text{for}\ 1\leq t\leq 3,\ \  \psi|_{\overline{P}^{S,c}_{e,5_t}(n)}:=\rho_t.
\end{equation*}

\textbf{Case 1:} Let $\lambda=(\lambda_{e_{1}},\lambda_{e_{2}},\dots,\lambda_{e_{r+1}})\cup(\lambda_{o_{1}},\lambda_{o_{2}},\dots,\lambda_{o_r}) \in \overline{P}^{S,c}_{e,5}(n)$. For $1\leq t\leq3$,
\begin{equation*}
\begin{split}
\rho_t(\lambda):=\lambda_e \cup (\lambda_{o_1}-2,\lambda_{o_{2}}-2,\dots,\lambda_{o_{r-t}}-2,1,1,\dots,1)\ \ \text{with}\ \mult_{\pi}(1)=2r-2t+\sum_{k=0}^{t-1}\lambda_{o_{r-k}}\\
\end{split}
\end{equation*}
and consequently for $1\leq t\leq 2$,
\begin{equation*}
\begin{split}
\rho_t(\overline{P}^{S,c}_{e,5_t}(n))&:= \overline{P}^{S,c}_{o,5_t}(n)\\
&:=\bigl\{\pi \in P^S_o(n): \mult_{\pi}(1)=2r-2t+\sum_{k=0}^{t-1}\lambda_{o_{r-k}},\ell_e(\pi)= r+1\ \text{and}\ \ell_{o,1}(\pi)=r-t\bigr\},
\end{split}
\end{equation*}
whereas for $t=3$,
\begin{equation*}
\begin{split}
\rho_3(\overline{P}^{S,c}_{e,5_3}(n))&:= \overline{P}^{S,c}_{o,5_3}(n)\\
&:=\bigl\{\pi \in P^S_o(n): \mult_{\pi}(1)=2r-2t+\sum_{k=0}^{2}\lambda_{o_{r-k}},\ell_e(\pi)= r+1\ \text{and}\ \ell_{o,1}(\pi)\leq r-3\bigr\}.
\end{split}
\end{equation*}
So for each $t$ we have, $$2r-2t+\sum_{k=0}^{t-1}\lambda_{o_{r-k}}>\dfrac{\ell(\pi)}{2}+1.$$
For example,
\begin{equation*}
\begin{split}
(10,9,7,7,6,4,4) & \xrightarrow{\rho_1} (10,7,6,5,4,4,1,1,1,1,1,1,1,1,1,1,1),\\
(10,9,7,6,4,4,3) & \xrightarrow{\rho_2} (10,7,6,4,4,1,1,1,1,1,1,1,1,1,1,1,1),\\
(10,9,6,4,4,3,3) & \xrightarrow{\rho_3} (10,6,4,4,1,1,1,1,1,1,1,1,1,1,1,1,1,1,1),\\
(10,9,8,7,6,6,4,4,4,3,3,3,3) & \xrightarrow{\rho_3} (10,8,7,6,6,5,4,4,4,1,1,1,1,1,1,1,1,1,1,1,1,1,1,1,1).
\end{split}
\end{equation*}

\textbf{Case 2:} For $n$ even, we have to consider $\overline{P}^{S,c}_{e,5_{4,2}}(n) \subset \overline{P}^{S,c}_{e,5_4}(n)$ as any partition $\lambda \in \overline{P}^{S,c}_{e,5_{4,2}}(n)$ is of the form $(\lambda_{e_1},\lambda_{e_2},\lambda_{e_3})\cup(3,3)$. Define $\rho_{4}(\lambda):=\pi=(\lambda_{e_{1}}+3,\lambda_{e_{2}}+3,\lambda_{e_{3}})$ so that $\ell_o(\pi)-\ell_{e}(\pi)=1$ and $a(\pi)$ is odd. One can observe that $\rho_4(\overline{P}^{S,c}_{e,5_{4,2}}(n))$ is disjoint from the sets which have already a preimage either by the map $f,\hat{f},\phi$ or by $\{\rho_t\}_{1\leq t\leq 3}$.

Similarly, if $n$ is odd, then we consider $\overline{P}^{S,c}_{e,5_{4,1}}(n) \subset \overline{P}^{S,c}_{e,5_4}(n)$. We compare the partitions $\lambda=(\lambda_{e_1},\lambda_{e_{2}})\cup(3) \in \overline{P}^{S,c}_{e,5_{4,1}}(n)$ with the partitions of the form $((\pi_{o_1},\pi_{o_2})\cup(3))\cup(4,4) \in P^S_o(n)$ where $\pi_{o_1}\geq 5$, and $((\pi_{o_1},\pi_{o_2})\cup(3))\cup(6,4)\in P^S_o(n)$ where $\pi_{o_1}\geq 7$, for all $n>23$. The total number of such partitions is
$$\bigg(\bigg\lfloor\frac{n-9}{4}\bigg\rfloor-1\bigg) + \bigg(\bigg\lfloor\frac{n-11}{4}\bigg\rfloor-1\bigg),$$
and
$$\Bigl|\overline{P}^{S,c}_{e,5_{4,1}}(n)\Bigr|=\bigg\lfloor\frac{n-3}{4}\bigg\rfloor-1.$$
Since $n$ is odd, it is of the form $4s+1$ or $4s+3$. For any form of $n$ we get, 
$$\bigg(\bigg\lfloor\frac{n-9}{4}\bigg\rfloor-1\bigg) + \bigg(\bigg\lfloor\frac{n-11}{4}\bigg\rfloor-1\bigg)>\bigg\lfloor\frac{n-3}{4}\bigg\rfloor-1.$$

We note that the set $\overline{P}^{S,c}_{e,5_{4,1}}(n)\neq \emptyset$ for $n=11,13,\dots,23$. In these cases $\lambda=(\lambda_{e_1},\lambda_{e_{2}})\cup(3)$ maps to $\rho_5(\lambda):=\pi=(\lambda_{e_2}-1,\lambda_{e_2}-1,1)\cup(\lambda_{e_1}-\lambda_{e_2}+4)$. $\ell_o(\pi)-\ell_e(\pi)=2$ and the two greatest odd parts are equal. So it was not an image in the above cases. For example, $(14,6,3)$ maps to $(12,5,5,1)$.

Now, we observe that $P_o^S(n)$ has more classes of partitions for all $n$, which are not mapped yet. For example, the partition where all parts are 1. Hence $p_o^{\{2\}}(n)>p_e^{\{2\}}(n)$ for all $n\geq1$.

This completes the proof. \qed

\emph{Proof of Theorem \ref{thm:missingset2}}: 
For $\lambda \in p^S_o(n)$, it is immediate that $\lambda_{e_i}\geq 4$, where $\lambda_e=(\lambda^{k_1}_{e_{1}}\lambda^{k_2}_{e_{2}}\dots\lambda^{k_s}_{e_{s}})$ is the even component of $\lambda$. So we can apply the injective map $f$ (cf. Section \ref{sec:proof-kkl}) on $G_e(n)\subseteq P_e(n)$. Consequently for $\lambda=(\lambda^{m_1}_{1}\lambda^{m_2}_{2}\dots\lambda^{m_k}_{k})\in G^0_e(n)$, divide the subset $G^0_e(n)$ into two disjoint classes given as follows
\begin{equation*}
\begin{split}
G^0_e(n) &:= G^0_{e,1}(n) \cup G^0_{e,2}(n),\\
\text{with}\ \hspace{2 cm} & \\
G^0_{e,1}(n) &= \{\lambda \in G^0_e(n): \lambda_{3}>6\}\\ 
\text{and}\ \ G^0_{e,2}(n) &= \{\lambda \in G^0_e(n): \lambda_{3}\in T\},\ \text{where}\  T=\{3,4,5,6\}.
\end{split}
\end{equation*} 
By a similar method to the proof of Theorem \ref{thm:reverse}, this result can be proved. We leave this as an exercise to the readers. \qed

\section{Concluding Remarks}\label{sec:final}

In continuation with the study on parity of parts, as discussed in the previous section (cf. Section \ref{sec:missingset}), we conclude this paper by presenting a somewhat more general discussion on the prospect towards further study.

First, by allowing only distinct partitions in $Q_e(n)$ and $Q_o(n)$ into Theorem \ref{thm:reverse}, we propose the following problem.
\begin{problem}
	For all $m>6$ we have $$dq_o(2m) > dq_e(2m),$$ and $$dq_o(2m+1) < dq_e(2m+1).$$
\end{problem}
\begin{remark}
	In fact, $dq_o(2m+1)<dq_e(2m+1)$ for $m=4,5$ as well.
\end{remark}
Whereas, Theorems \ref{thm:missingset1} and \ref{thm:missingset2} suggest a more general family of inequalities in the following sense. 
\begin{problem}
 For all $k>2$ we have $p_o^{\{k\}}(n)>p_e^{\{k\}}(n)$ and $p_e^{\{1,k\}}(n)>p_o^{\{1,k\}}(n)$, for all $n>N(k)$, for some constant $N(k)$, depending on $k$.
 
Moreover, it would be worthwhile to understand the threshold $N(k)$ asymptotically.
\end{problem}

\appendix

\section{Proofs of Lemmas \ref{lem1} and \ref{lem2}}\label{sec:appendix}

\emph{Proof of Lemma \ref{lem1}}:
	Take $n = 2m$ with $m \in \mathbb{Z}_{\geq 7}$. So, we have to show
	\begin{equation}\label{fund1}
	\sum_{k=1}^{m-3}\Bigl\lfloor\frac{m-k-1}{2}\Bigr\rfloor > 1 + \sum_{k=1}^{\lfloor \frac{m-1}{3}\rfloor}\Bigl\lfloor \frac{m-3k+1}{2}\Bigr\rfloor + \sum_{k=1}^{\lfloor \frac{m-3}{3}\rfloor}\Bigl\lfloor \frac{m-3k-1}{2}\Bigr\rfloor.
	\end{equation}
	Note that
	\begin{equation*}
	\begin{split}
	\sum_{k=1}^{m-3}\Bigl\lfloor\frac{m-k-1}{2}\Bigr\rfloor = \sum_{k=0}^{m-4}\Bigl\lfloor\frac{k+2}{2}\Bigr\rfloor & >\sum_{k=0}^{m-4}\Bigl(\frac{k+2}{2}-1\Bigr) \ \Bigl(\text{since,} \ \lfloor x\rfloor > x-1\Bigr)\\
	& = \frac{(m-4)(m-3)}{4} = \frac{m^2-7m+12}{4},\\
	\end{split}
	\end{equation*}
	and
	\begin{equation*}
	\begin{split}
	&1 + \sum_{k=1}^{\lfloor \frac{m-1}{3}\rfloor}\Bigl\lfloor \frac{m-3k+1}{2}\Bigr\rfloor + \sum_{k=1}^{\lfloor \frac{m-3}{3}\rfloor}\Bigl\lfloor \frac{m-3k-1}{2}\Bigr\rfloor \\
	& < 1 + \sum_{k=1}^{\lfloor \frac{m-1}{3}\rfloor} \frac{m-3k+1}{2} + \sum_{k=1}^{\lfloor \frac{m-3}{3}\rfloor} \frac{m-3k-1}{2} \ \Bigl(\text{since,} \ \lfloor x\rfloor < x\Bigr)\\
	&= 1+\frac{m+1}{2}\Bigl\lfloor\frac{m-1}{3}\Bigr\rfloor- \frac{3}{4}\Bigl\lfloor\frac{m-1}{3}\Bigr\rfloor\Bigl(\Bigl\lfloor\frac{m-1}{3}\Bigr\rfloor+1\Bigr)+\frac{m-1}{2}\Bigl\lfloor\frac{m-3}{3}\Bigr\rfloor\\ & \qquad -\frac{3}{4}\Bigl\lfloor\frac{m-3}{3}\Bigr\rfloor\Bigl(\Bigl\lfloor\frac{m-3}{3}\Bigr\rfloor+1\Bigr)\\
	& < 1+\frac{m+1}{2}\frac{m-1}{3}- \frac{3}{4}\Bigl(\frac{m-1}{3}-1\Bigr)\Bigl(\frac{m-1}{3}\Bigr)+\frac{m-1}{2}\frac{m-3}{3}- \frac{3}{4}\Bigl(\frac{m-3}{3}-1\Bigr)\Bigl(\frac{m-3}{3}\Bigr)\\
	& = \frac{m^2+3m-3}{6}.\\
	\end{split}
	\end{equation*} 
	Now, $\dfrac{m^2-7m+12}{4}> \dfrac{m^2+3m-3}{6}$ for $m \in \mathbb{Z}_{\geq 26}$. We finish the proof by checking the inequality \eqref{fund1} for $7 \leq m \leq 25$ numerically in Mathematica.\qed
\\

\emph{Proof of Lemma \ref{lem2}}:
	The proof is exactly similar to the proof of Lemma \ref{lem1}.\qed

\section*{Acknowledgements}
Authors would like to thank the anonymous referees for helping to improve the first version of the paper. The first author would like to acknowledge that the research was funded by the Austrian Science Fund (FWF): W1214-N15, project DK6. The last author is partially supported by the Leverhulme Trust Research Project Grant RPG-2019-083 and thanks Prof. Michael Schlosser (Vienna) for suggesting to look at parity biases in partitions with restrictions on the parts.

\end{document}